\documentclass[12pt,a4paper]{article}
\usepackage{amsmath}
\usepackage{amssymb}
\usepackage{amsthm}
\usepackage{amsfonts}
\usepackage[all]{xy}
\usepackage{cite}
\usepackage{hyperref}
\usepackage{color}
\usepackage{tikz}

\newtheorem{theorem}{Theorem}
\newtheorem{proposition}[theorem]{Proposition}

\newtheorem{lemma}[theorem]{Lemma}
\newtheorem{corollary}[theorem]{Corollary}

\newtheorem{claim}{Claim}

\newtheorem{observation}[theorem]{Observation}

\def\tobedone{\textcolor{red}{TO BE DONE.}}

\def\diam{\operatorname{diam}}

\renewenvironment{proof}{
\par
\noindent {\bf Proof.}\rm}{\mbox{}\hfill$\square$\par\vskip 3mm}
\newcommand\EFFACE[1]{}

\def\diam{{\rm diam}}

\newcommand\LIGNE[4]{
\draw[thick] (#1,#2) to (#3,#4);
}
\newcommand\POINTILLE[4]{
\draw[thick,dotted] (#1,#2) to (#3,#4);
}

\newcommand\bSOM[4]{
   \node[scale=0.7,draw,circle,fill=black] at (#1,#2){};
   \node[above] at (#1,#2+0.2){#3};
   \node[below] at (#1,#2-0.2){#4};
}
\newcommand\gSOM[4]{
   \node[scale=0.7,draw,circle,fill=lightgray] at (#1,#2){};
   \node[above] at (#1,#2+0.2){#3};
   \node[below] at (#1,#2-0.2){#4};
}
\newcommand\SOM[4]{
   \node[scale=0.7,draw,circle,fill=white] at (#1,#2){};
   \node[above] at (#1,#2+0.2){#3};
   \node[below] at (#1,#2-0.2){#4};
}

\begin{document}

\title{Broadcast independence number \\ of oriented circulant graphs}

\author{Abdelamin LAOUAR~\thanks{Faculty of Mathematics, Laboratory L'IFORCE, University of Sciences and Technology
Houari Boumediene (USTHB), B.P.~32 El-Alia, Bab-Ezzouar, 16111 Algiers, Algeria.}
\and Isma BOUCHEMAKH~\footnotemark[1]
\and \'Eric SOPENA~\thanks{Univ. Bordeaux, CNRS, Bordeaux INP, LaBRI, UMR5800, F-33400 Talence, France.}
}
\maketitle

\begin{abstract}
In 2001, D. Erwin \cite{Erw01} introduced in his Ph.D. dissertation  the notion of broadcast independence in unoriented graphs. Since then, some results but not many, are published on this notion, including research work on the broadcast independence number  of  unoriented circulant graphs \cite{LBS23}. In this paper, we are focused in the same parameter but of the class of oriented  circulant graphs.\\
An independent broadcast on an oriented graph $\overrightarrow{G}$ is a function  $f: V\longrightarrow \{0,\ldots,\diam(\overrightarrow{G})\}$ such that
$(i)$ $f(v)\leq e(v)$ for every vertex $v\in V(\overrightarrow{G})$, where $\diam(\overrightarrow{G})$ denotes the diameter of $\overrightarrow{G}$ and $e(v)$ the eccentricity of vertex $v$, and $(ii)$ $d_{\overrightarrow{G}}(u,v) > f(u)$ for every distinct vertices $u$, $v$ with  $f(u)$, $f(v)>0$, where $d_{\overrightarrow{G}}(u,v)$ denotes the length of a shortest oriented path from $u$ to $v$. The broadcast independence number $\beta_b(\overrightarrow{G})$ of $\overrightarrow{G}$ is then the maximum value of $\sum_{v \in V} f(v)$, taken over all independent broadcasts on $\overrightarrow{G}$.\\ The goal of this paper is to study the properties of independent broadcasts of oriented circulant graphs $\overrightarrow{C}(n;1,a)$, for any integers $n$ and $a$ with $n>|a|\geq 1$ and $a \notin \{1,n-1\}$. Then, we give some bounds and some exact values for the number $\beta_b(\overrightarrow{C}(n;1,a))$.
\end{abstract}

\medskip


\section{Introduction}\label{sec:intro}
All the graphs we considered are simple. For such a graph $G$,  we denote by $V(G)$ and $E(G)$ its set of vertices and its set of edges, respectively.
For a nontrivial undirected graph $G$, the \emph{distance} from a vertex $u$ to a vertex $v$ in $G$, denoted $d_G(u,v)$, or simply $d(u,v)$ when
$G$ is clear from the context, is the length (number of edges) of a shortest $u-v$ path (by a $u-v$ path in $G$, we mean a path in $G$ whose end-vertices are the
vertices $u$ and $v$). The \emph{eccentricity} of a vertex $v$ in $G$, denoted $e_G(v)$, is the maximum distance from $v$ to any other vertex of $G$.
The  maximum eccentricity in $G$ is its \emph{diameter}, denoted $\diam(G)$.\\
A function $f: V(G) \rightarrow \{0,\ldots,\diam(G)\}$ is a \emph{broadcast on $G$} if $f(v) \leq e_G(v)$ for every vertex $v \in V$. For each vertex $v$, $f(v)$ is the \emph{$f$-value} of $v$, or the \emph{broadcast value} of $v$ if $f$ is clear from the context.
Given such a broadcast $f$, an \emph{$f$-broadcast vertex} is a vertex $v$ for which $f(v) > 0$. The set of all $f$-broadcast vertices is denoted $V^+_f(G)$. If $v$ is a broadcast vertex and $u$ a vertex such that $d(u,v) \leq f(v)$, then the vertex $v$ \emph{$f$-dominates} the vertex $u$. The \emph{cost} of a broadcast $f$ on $G$ is the value $\sigma(f) =\sum_{v\in V_f^+} f(v)$. \\
A broadcast $f$ is \emph{independent} if no broadcast vertex $f$-dominates another broadcast vertex, or, equivalently, if $d(u,v) > \max\{f(u), f(v)\}$ for every two distinct broadcast vertices $u$ and $v$. The maximum cost of an independent broadcast on $G$ is the \emph{broadcast independence number} of $G$, denoted $\beta_b(G)$.
An independent broadcast with cost $\beta_b(G)$ is referred to as a \emph{$\beta_b$-broadcast}.
The notion of independence broadcast extends to digraphs in a natural way:
Let $\overrightarrow{D}$ be a digraph, with vertex set $V(\overrightarrow{D})$ and arc set $E(\overrightarrow{D})$.  A \emph{directed path} of length $k$ in $\overrightarrow{D}$ is a sequence $u_0\ldots u_k$ of vertices of $V(\overrightarrow{D})$ such that for every $i$, $0\leq  i \leq  k - 1$, $u_iu_{i+1}$ is an arc in $E(\overrightarrow{D})$. Such a path is denoted $u_0-u_k$ path. The \emph{weak directed distance} between two vertices $u$ and $v$ in $\overrightarrow{D}$, denoted $d_{\overrightarrow{D}}(u, v)$, is the shortest length (number of arcs) of a directed $u-v$ path or $v-u$ path in $\overrightarrow{D}$. A digraph $\overrightarrow{O}$ with no pair of opposite arcs, that is $uv\in  E(\overrightarrow{O})$ implies $vu \notin E(\overrightarrow{O})$, is called an \emph{oriented graph}. If $G$ is an undirected graph, an \emph{orientation} of $G$
is any oriented graph $\overrightarrow{G}$ obtained by giving to each edge of $G$ one of its two possible orientations.\\

In this paper, we study the broadcast independence number of oriented circulant graphs. We first define this class of graphs in the unoriented case. For every integer $n\ge 3$, and every sequence of integers $a_1,\dots,a_k$, $k\ge 1$, satisfying $1 \leq a_1 \leq \dots \leq a_k\leq \left\lfloor \frac{n}{2}\right\rfloor$, the \emph{circulant graph} $G = C(n;a_1,\dots,a_k)$ is the graph defined by
$$V(G) = \{v_0,v_1,\dots,v_{n-1}\}\ \mbox{and }
E(G) = \left\{v_iv_{i+a_j}\ |\ a_j\in \{a_1,\dots,a_k\} \right\}$$
(subscripts are taken modulo $n$).

We now define this class of circulant graph in oriented case:  For any integers $n$ and $a_1,\ldots,a_k$, with $a_i+ a_j \neq n$, for every $i,j=1,\ldots,k$ the \emph{oriented circulant graph}  $\overrightarrow{G}= \overrightarrow{C}(n;a_1,\dots, a_k)$ is the graph defined by $$V(\overrightarrow{G}) = \{v_0,v_1,\dots,v_{n-1}\}\ \mbox{and }
 E(\overrightarrow{G}) = \left\{v_iv_{i+a_j}\ |\ a_j\in \{a_1,\dots,a_k\} \right\}$$
(subscripts are taken modulo $n$).

Broadcast independence was introduced by Erwin~\cite{Erw01} in his Ph.D. dissertation, using the term cost independence. He also discussed several other types of broadcast parameters and gave relationships between them. Most of the corresponding results are published in~\cite{DEHHH,Erw04}. Since then, some results but not many, are published on the broadcast independence number   (see \cite{ABS18},\cite{ABS19}, \cite{BZ12,BBS19,DEHHH},\cite{LBS23}), on the algorithmic complexity of broadcast independence \cite{BR1_2018} and on the links between girth, minimum degree, independence number and broadcast independence number \cite{BR2_2018,BR3_2018}.\\
In \cite{LBS23}, we  prove that every (unoriented) circulant graph of the form $C(n; 1, a)$, $3\leq a \leq \lfloor\frac{n}{2}\rfloor$, admits an optimal 2-bounded independent broadcast, that is, an independent broadcast $f$ satisfying $f(v)\leq 2$ for every vertex $v$, except when $n = 2a+1$,  or $n=2a$ and $a$ is even. We determine the broadcast independence number of various classes of such circulant graphs, and prove that, for most of these classes, the equality $\beta_b(C(n; 1, a))=\alpha(C(n; 1, a))$ holds, where $\alpha(C(n; 1, a))$ denotes the independence number of $C(n; 1, a)$.

In this study, we  focus on oriented circulant graph $\overrightarrow{C}(n;1,a)$, with  $n>|a|\geq 1$ and $a \notin \{1,n-1\}$. We define the distance $d_{\overrightarrow{G}}(u, v)$, or simply $d(u,v)$ when
$\overrightarrow{G}$ is clear from the context, between two vertices $u$ and $v$ in  $\overrightarrow{C}(n;1,a)$, as the shortest length (number of arcs) of a directed path in $\overrightarrow{G}$ going either from $u$ to $v$ (without take into account of shortest length  of a directed path going either from $v$ to $u$ and this restriction has no impact on our results).

If $u$ is a broadcast vertex and $v$ a vertex such that $d_{\overrightarrow{G}}(u,v) \leq f(u)$, then the vertex $u$ \emph{$f$-dominates} the vertex $v$. The \emph{cost} of a broadcast $f$ on $\overrightarrow{G}$ is the value $\sigma(f) =\sum_{v\in V_f^+} f(v)$. \\

A function $f$ is an \emph{independent broadcast} on $\overrightarrow{G}=\overrightarrow{C}(n;1,a)$  if  $f(v)\leq e_{\overrightarrow{G}}(v)$ for every vertex $v\in V(\overrightarrow{G})$ and  $d_{\overrightarrow{G}}(u, v) > f(u)$ for every distinct vertices $u$, $v$ with  $f(u)$, $f(v)>0$. The maximum cost of an independent broadcast on $\overrightarrow{G}$ is the \emph{broadcast independence number} of $\overrightarrow{G}$, denoted $\beta_b(\overrightarrow{G})$. An independent broadcast with cost $\beta_b(\overrightarrow{G})$ is referred to as a \emph{$\beta_b$-broadcast}.\\

\begin{figure}[h]
\begin{center}
\begin{tikzpicture}[scale=0.6]

\draw (0, 0) circle (4cm);

\node[rotate=-22] at (1,3.87) {\textbf{{\large $>$}}};
\node[rotate=-52] at (2.84,2.84) {\textbf{{\large $>$}}};
\node[rotate=-55] at (1.6,1.3) {\textbf{{\large $>$}}};
\node[rotate=-180] at (0,-2) {\textbf{{\large $>$}}};
\node[rotate=-140] at (2.85,-2.8) {\textbf{{\large $>$}}};

\LIGNE {0}{4}{3.5}{-2}
\LIGNE {3.5}{-2}{-3.5}{-2}

\POINTILLE{0}{4}{-3.5}{-2}

\POINTILLE {0}{4}{3.5}{-2}
\POINTILLE{0}{4}{-3.5}{-2}
\POINTILLE {3.5}{-2}{-3.5}{-2}

\POINTILLE {2}{3.5}{2}{-3.5}
\POINTILLE{2}{-3.5}{-4}{0}
\POINTILLE{2}{3.5}{-4}{0}
\POINTILLE {3.5}{2}{0}{-4}
\POINTILLE {-3.5}{2}{0}{-4}
\POINTILLE {3.5}{2}{-3.5}{2}
\POINTILLE {-2}{3.5}{4}{0}
\POINTILLE {-2}{3.5}{-2}{-3.5}
\POINTILLE {-2}{-3.5}{4}{0}
 \bSOM{0}{4}{$2$}{}
  \gSOM{2}{3.5}{\quad $0$}{}
 \gSOM{3.5}{2}{\quad $0$}{}
  \bSOM{4}{0}{ }{}
  \node[color=black,thick]at(4.7,0){$2$};
   \gSOM{3.5}{-2}{}{\quad $0$}
   \gSOM{2}{-3.5}{}{\quad $0$}
  \bSOM{0}{-4}{}{$2$}
  \gSOM{-2}{-3.5}{}{ $0$ \quad }
   \gSOM{-3.5}{-2}{}{}
    \node[color=black,thick]at(-3.8,-2.5){$0$};
    \bSOM{-4}{0}{ }{}
 \node[color=black,thick]at(-4.7,0){$2$};
 \gSOM{-3.5}{2}{}{  }
  \node[color=black,thick]at(-3.9,2.5){$0$};
   \gSOM{-2}{3.5}{$0$ \quad }{}
\node[below] at (0,-6){$(a): \sigma(f_1)= 8 $};
 \draw (12, 0) circle (4cm);
 \node[rotate=-22] at (13,3.87) {\textbf{{\large $>$}}};
\node[rotate=-52] at (14.84,2.84) {\textbf{{\large $>$}}};
\node[rotate=-35] at (13.8 ,3) {\textbf{{\large $>$}}};
 \node[rotate=-90] at (15.5 ,0) {\textbf{{\large $>$}}};
\node[rotate=-145] at (13.8 ,-3) {\textbf{{\large $>$}}};
\node[rotate=145] at (10.2 ,-3) {\textbf{{\large $>$}}};
\node[rotate=90] at (8.5 ,0) {\textbf{{\large $>$}}};
 \node[rotate=-90] at (15.5 ,0) {\textbf{{\large $>$}}};
 \node[rotate=-80] at (15.9,1) {\textbf{{\large $>$}}};
 \node[rotate=-110] at (15.85,-1) {\textbf{{\large $>$}}};
 \node[rotate=-140] at (14.85,-2.8) {\textbf{{\large $>$}}};
 \node[rotate=-165] at (13.2,-3.8) {\textbf{{\large $>$}}};

\LIGNE {12}{4}{15.5}{2}
\LIGNE {15.5}{-2}{15.5}{2}
\LIGNE {15.5}{-2}{12}{-4}
\LIGNE {8.5}{-2}{12}{-4}
\LIGNE {8.5}{-2}{8.5}{ 2}
\POINTILLE {14}{3.5}{16}{0}
\POINTILLE {14}{3.5}{10}{3.5}
\POINTILLE {15.5}{-2}{15.5}{ 2}
\POINTILLE {14}{-3.5}{10}{-3.5}
\POINTILLE {14}{-3.5}{16}{0}
\POINTILLE {8.5}{-2}{8.5}{ 2}
\POINTILLE {8.5}{-2}{12}{-4}
\POINTILLE {10}{-3.5}{8}{0}
\POINTILLE {10}{3.5}{8}{0}
 \bSOM{12}{4}{$6$}{}
  \gSOM{14}{3.5}{\quad $0$}{}
 \gSOM{15.5}{2}{\quad $0$}{}

  \gSOM{16}{0}{ }{}
  \node[color=black,thick]at(16.7,0){$0$};
   \gSOM{15.5}{-2}{}{\quad $0$}
   \gSOM{14}{-3.5}{}{\quad $0$}
  \gSOM{12}{-4}{}{$0$}
  \gSOM{10}{-3.5}{}{ $0$ \quad  \quad }
   \gSOM{8.5}{-2}{}{}
    \node[color=black,thick]at(8,-2.5){$0$};
    \gSOM{8}{0}{ }{}
 \node[color=black,thick]at(7.3,0){$0$};
 \gSOM{8.5}{2}{}{  }
  \node[color=black,thick]at(8.1,2.5){$0$};
   \gSOM{10}{3.5}{$0$ \quad }{}
 \node[below] at (12,-6){$(b): \sigma(f_2)=6$};
\end{tikzpicture}
\caption{\label{fig:broadcast} Two independent broadcasts $f_1$ on $\overrightarrow{C}(12;1,4)$ and $f_2$ on $\overrightarrow{C}(12;1,2)$.}
\end{center}

\end{figure}
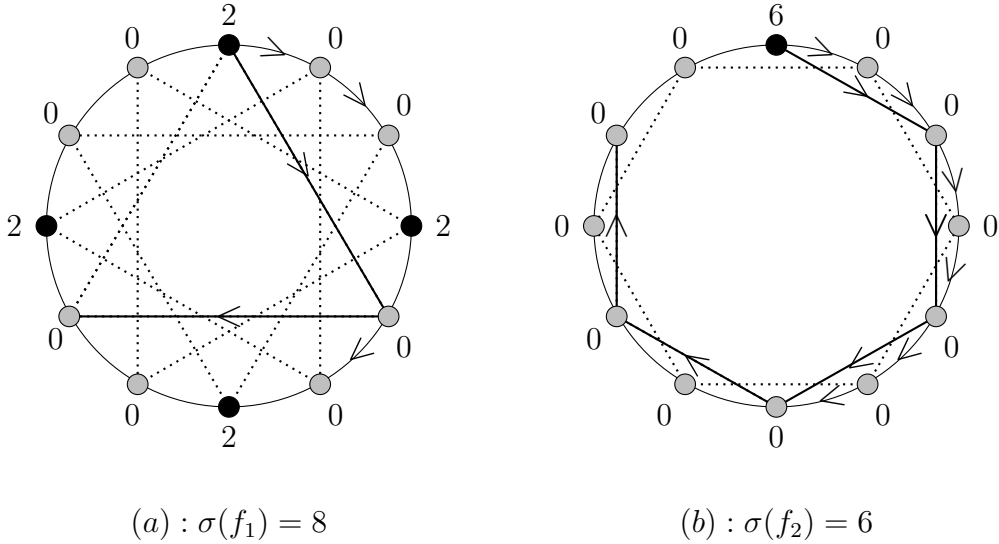

Figures~\ref{fig:broadcast}(a) and ~\ref{fig:broadcast}(b)  illustrate two independent broadcasts  $f_1$ on $\overrightarrow{C}(12;1,4)$ with cost $\sigma(f_1) = 8$ and $f_2$ on $\overrightarrow{C}(12;1,2)$  with cost $\sigma(f_2) = 6$.

By definition, if $\overrightarrow{G}$ is any orientation of an undirected graph $G$ then, for any two vertices $u$ and $v$ in $G$, $ d_G(u,v) \le d_{\overrightarrow{G}}(u,v) $. Consequently, every independent broadcast of $G$ is a independent broadcast of $\overrightarrow{G}$. This leads to the following bound. Moreover, this bound is achieved for some classes of oriented circulant graphs, which we introduce in the next section.
\begin{proposition}\label{prop:lower bound 1}
If $n$ and $a$ are two integers such that $n > 3$, $|a|\geq 1$, and $a \notin\{1,n-1\}$, then
$$ \beta_b(\overrightarrow{G})  \geq \max \{ \beta_b (G), \diam (\overrightarrow{G})\}.$$
\end{proposition}

Our paper is organized as follows.
After general introduction, we introduce in Section 2 few isomorphisms between circulant graphs, allowing us to reduce our workload a little bit, and we give exact values of $\beta_b(\overrightarrow{C}(n;1,a))$ for $n=2a$, $n=2a-1$, and $n=3a-1$ with $a\geq 4$, and for every $n$ with  $a=2$ or $a=3$. In Section 3, we identify some classes of oriented circulant graphs $\overrightarrow{C}(n;1,a)$ for which there exists an $\ell$-bounded independent $\beta_b$-broadcast. Then, we give exact values of $\beta_b(\overrightarrow{C}(n;1,a))$ for $n=k(a-1)$ and $k\geq 3$, for $n=qa$ with $q=k(a-1)$ and $k\geq 1$, or $a=kq+1$ and $k\geq 2$ and for $n=qa+a-1$ and $k\geq 1$.

\section{Preliminary results}\label{sec:preliminaries}

For any integers $n$ and $a$ with $n > |a|  \geq 1$, $a \notin \{1,n-1\}$, let us point out that  the graphs $\overrightarrow{C}(n;1,a)$ and $\overrightarrow{C}(n;-1,-a)$ are isomorphic, as well as the graphs $\overrightarrow{C}(n;-1,a)$ and $\overrightarrow{C}(n;1,-a)$, leading immediately to the relations $\beta_b (\overrightarrow{C}(n;1,a))= \beta_b ( \overrightarrow{C}(n;-1,-a))$ and $\beta_b ( \overrightarrow{C}(n;-1,a) )= \beta_b ( \overrightarrow{C}(n;1,-a))$. Based on this observation, we can restrict our study to only on oriented circulant graphs $\overrightarrow{C}(n;1,\pm a)$. Another isomorphic reduces a little more this study, as showed in the following observation
\begin{observation}\label{isomorphisme}
If $n$ and $a$ are two integers such that, $2 \le a \le n-2 $, then
$$\overrightarrow{C}(n;1,a) \simeq \overrightarrow{C}(n;1,-(n-a)).$$
\end{observation}
Thanks to Observation~\ref{isomorphisme}, we can streamline our focus to delve exclusively into the properties and characteristics of the graph $\overrightarrow{C}(n;1,a)$ with $a\geq 2.$

\begin{theorem}\label{th:C(2a;1,a)}
For every integer $a\ge 2$,
$$\beta_b(\overrightarrow{C}(2a;1,a)) = \beta_b(\overrightarrow{C}(2a;1,-a)) =  \diam (\overrightarrow{C}(2a;1,a)) = a.$$
\end{theorem}
\begin{proof}
Let $f$ be an independent broadcast on $\overrightarrow{C}(2a;1,a)$ and let
$V_f^{1}=\{v_i\in V_f\ |\ f(v_i) = 1\}$ and $V_f^{\geq 2}=\{v_i\in V_f\ |\ f(v_i) \geq 2 \}$ a partition of $V_f^+$.\\
 For each vertex $v_i\in V_f^1$ and $v_j\in V_f^{\geq 2}$,  we set
\begin{center}
$ A_f^i = \{v_{i},v_{i+1}\} $ and
$ B_f^j = \{v_{j}, v_{j+1}, \ldots v_{j+f(v_j)}\} \cup \{v_{j+a+1}, v_{j+a+2}, \ldots v_{j+a+f(v_j)-1}\}$.
\end{center}
The definition of the set $B_f^j$ is illustrated in Figure~\ref{fig:Bj} with $n=16$ and $a = 8$. \\
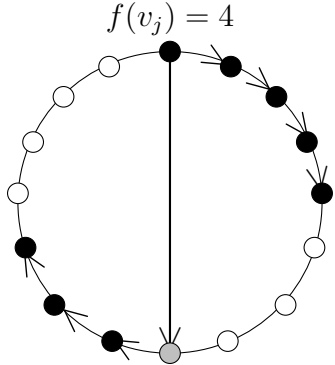
\begin{figure}[h]
\begin{center}
\begin{tikzpicture}[scale=0.4]

\draw(0,0)circle(5) ;
\LIGNE{0}{5}{0}{-5}
\node[rotate=-90] at (0,-4.5) {\textbf{{\large $>$}}};
   \bSOM{0}{5}{$f(v_j)=4$}{}
   \bSOM{2}{4.5}{}{}
   \bSOM{3.5}{3.5}{}{}
   \bSOM{4.5}{2}{}{}
   \bSOM{5}{0.3}{}{}
   \SOM{4.75}{-1.5}{}{} \node[below] at (5.5,-0.8) { };
   \SOM{3.8}{-3.4}{}{}
   \SOM{1.9}{-4.6}{}{}
   \gSOM{0}{-5}{}{}
  \bSOM{-1.9}{-4.6}{}{}
  \bSOM{-3.8}{-3.4}{}{}
   \bSOM{-4.75}{-1.5}{}{}
   \SOM{-5}{0.3}{}{}
    \SOM{-4.5}{2}{}{}
  \SOM{-3.5}{3.5}{}{}
     \SOM{-2}{4.5}{}{}
  \node[rotate=-30] at (1.5,4.75) {\textbf{{\large $>$}}};
  \node[rotate=-40] at (3,3.98) {\textbf{{\large $>$}}};
 \node[rotate=150] at (-1.5,-4.75) {\textbf{{\large $>$}}};
  \node[rotate=140] at (-3.2,-3.8) {\textbf{{\large $>$}}};

  \node[rotate=-54] at (4.3,2.5) {\textbf{{\large $>$}}};
  \node[rotate=-80] at (4.9,0.8) {\textbf{{\large $>$}}};
  \node[rotate=120] at (-4.6,-2) {\textbf{{\large $>$}}};

\tobedone
\end{tikzpicture}
\caption{\label{fig:Bj} The  set $B_f^j$ (the black vertices) with $n=16$ and $a=8$.}
\end{center}
\end{figure}

It is clearly seen that these sets are pairwise disjoint with sizes $|A_f^i|=2$ and $|B_f^j|=2 f(v_j)$. Then
$$\sum_{v_i\in V_f^1}|A_f^i|  + \sum_{v_j\in V_f^{\geq 2}}|B_f^j|
   =  2f(V_f^1)  + 2 f(V_f^{ \geq 2}) \le n,$$
which gives
$$ f(V_f^+) = f(V_f^1) + f(V_f^{\geq 2}) \le \frac{n}{2}=a.$$

Now, since in $\overrightarrow{C}(2a;1,a)$, $\diam (\overrightarrow{C}(2a;1,a))=e(v_0)=d(v_0,v_{n-1})= a$, we get, by Proposition~\ref{prop:lower bound 1},  $\beta_b(\overrightarrow{C}(2a;1,a)) \geq a$, which implies $\beta_b(\overrightarrow{C}(2a;1,a)) = a$. \\By Observation~\ref{isomorphisme}, we get $\beta_b(\overrightarrow{C}(2a;1,a)) = \beta_b(\overrightarrow{C}(2a;1,-a)) =  \diam (\overrightarrow{C}(2a;1,a)) = a.$ This completes the proof.
\end{proof}

\begin{theorem}\label{th:C(3a-1;1,a)}  
If $n$ and $a$ are two integers such that $n=3a-1$ and $ a \geq 4$, then
\begin{center}
$\beta_b(\overrightarrow{C}(n;1,a)) = \beta_b(\overrightarrow{C}(n;1,-(2a-1))) = \left\lfloor\frac{n}{2}\right\rfloor.$
\end{center}
\end{theorem}

\begin{proof}
Let $f$ be an independent $\beta_b$-broadcast on $\overrightarrow{C}(n;1,a)$.\\
For every $f$-broadcast vertex $v_i\in V_f^+$, let $L(v_i)$ be the set of vertices that are $f$-dominated only by $v_i$,that is
 $$ L(v_i) =  \left\{ v_{i}, v_{i+1}, \ldots ,v_{i+f(v_{i})} \right\}
 \cup \left\{ v_{i+a+1}, v_{i+a+2}, \ldots ,v_{i+a+f(v_{i})-1} \right\}\cup \left\{ v_{i+2a}, v_{i+2a+1}, \ldots ,v_{i+2a+f(v_{i})-2} \right\}.$$
Figure~\ref{fig:L(vi)} illustrates $L(v_i)$ on  oriented circulant graph  $\overrightarrow{C}(20;1,7)$ for a broadcast vertex $v_i$ with $f(v_i)=4$.

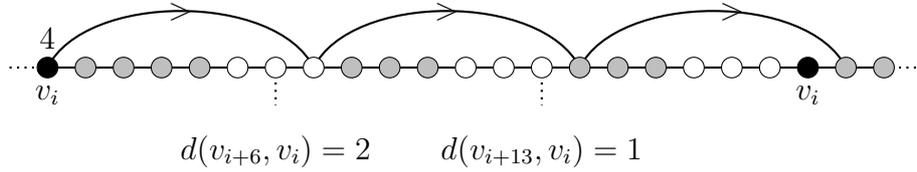
\begin{figure}[h]
\begin{center}
\begin{tikzpicture}[scale=0.5]
   \POINTILLE{3}{0}{4}{0}
   \LIGNE{4}{0}{26}{0}
   \POINTILLE{26}{0}{27}{0}
 \POINTILLE{10}{0}{10}{-1}
  \POINTILLE{17}{0}{17}{-1}
   \foreach \k in {4,11,18}
      \draw[thick] (\k,0) .. controls (\k+1,2) and (\k+6,2) .. (\k+7,0);
   \bSOM{4}{0}{4}{$v_i$}
   \gSOM{5}{0}{}{}
   \gSOM{6}{0}{}{}
   \gSOM{7}{0}{}{}
   \gSOM{8}{0}{}{}
   \SOM{9}{0}{}{}
   \SOM{10}{0}{}{}
   \SOM{11}{0}{}{}
   \gSOM{12}{0}{}{}
  \gSOM{13}{0}{}{}
  \gSOM{14}{0}{}{}
  \SOM{15}{0}{}{}
   \gSOM{26}{0}{}{}
   \gSOM{25}{0}{}{}
   \bSOM{24}{0}{}{$v_i$}
   \SOM{23}{0}{}{}
   \SOM{22}{0}{}{}
   \SOM{21}{0}{}{}
   \gSOM{20}{0}{}{}
   \gSOM{19}{0}{}{}
   \gSOM{18}{0}{}{}
   \SOM{17}{0}{}{}
   \SOM{16}{0}{}{}
\tobedone
 \node[below] at (7.5,2) {$>$};
 \node[below] at (14.5,2) {$>$};
 \node[below] at (22,2) {$>$};
\node[below] at (10,-1.5) {$d(v_{i+6},v_{i})=2$};
 \node[below] at (17,-1.5) {$d(v_{i+13}, v_{i})=1$};
\end{tikzpicture}
\caption{\label{fig:L(vi)}The set $L(v_i)$ (black vertex and grey vertices), with $a=7$ and $f(v_i)=4$.}
\end{center}
\end{figure}
 The sets $L(v_i)$ have the following two properties :
\begin{enumerate}
\item $L(v_i) \cap L(v_j) = \emptyset$ for every
two distinct vertices $v_i$ and $v_{j}$ in $V_f^+$, since otherwise we would have $d(v_{i},v_{j})\leq \max \{f(v_{i}),f(v_{j})\}$,
contradicting the fact that $f$ is an independent broadcast.
\item For every vertex $v_i\in V_f^+$, $|L(v_i)| = 3f(v_i)-1$.
\end{enumerate}
From these two properties we deduce that

$$\sum_{v_i\in V_f^+}|L(v_i)|  =  3f(V_f^+) - \left|V_f^+\right| \leq n,$$
which gives

$$f( V_f^+) \leq \left\lfloor   \frac{n + \left|V_f^+\right|}{3} \right \rfloor .$$

Finally, since $\left|V_f^{+}\right| \leq  \left\lfloor  \dfrac{n}{2} \right \rfloor $, we get

$$\sigma(f) = f( V_f^+) \leq \left\lfloor \frac{n }{2} \right \rfloor.$$

For the reverse inequality, we construct a mapping $g$ from $V(\overrightarrow{C}(n;1,a))$ to $\{0,1,2\}$
 with cost $ \left\lfloor  \dfrac{n}{2} \right \rfloor$. For this, we consider two cases, depending on the parity of  $a$.

\begin{enumerate}
\item If $a$ is odd, we let  $g(v_i)=1$ if  $i$ is even, and $g(v_i)=0$ otherwise. Since $a$ is odd, then $n$ is even and $g$ is clearly an independent broadcast on $\overrightarrow{C}(n;1,a)$ with cost $\sigma(g)=\frac{n}{2}$.
Hence $\beta_b(\overrightarrow{C}(n;1,a)) \ge \frac{n}{2}$.

\item If $a$ is even, we let
$$ g(v_i) = \left\{
  \begin{array}{ll}
     2 , & \mbox{if $i=0$,} \\ [1ex]
    1, & \mbox{if ($i$ is odd, and  $3\leq i \leq a-1 $ or $2a + 1 \leq i \leq 3a-3 $),} \\ [1ex]
    {} & \mbox{or if ($i$ is even, and  $a+2\leq i \leq 2a -2  $),} \\ [1ex]
    0, & \mbox{otherwise.}
  \end{array}
\right.$$
(See Figure~\ref{fig:mapping 1} for the case $n=17$ and $a=6$)

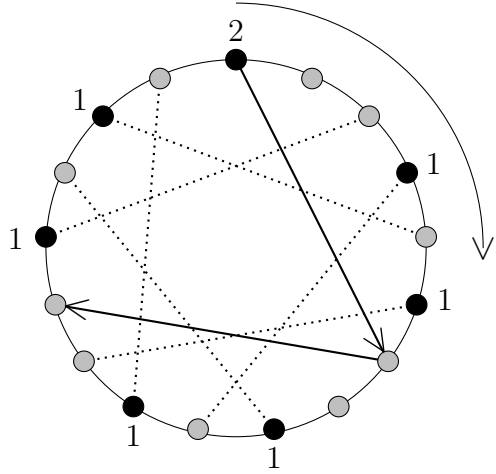
\begin{figure}[h]
\begin{center}
\begin{tikzpicture}[scale=0.5]
\draw(0,0)circle(5) ;
\draw (0,6.5) arc (90: 0 : 6.5cm);
   \node[rotate=-90] at (6.5,0) {\textbf{{\large $>$}}};
\LIGNE{0}{5}{4}{-3}
\LIGNE{4}{-3}{-4.75}{-1.5}
   \POINTILLE{4.5}{2}{-1}{-4.8}
   \POINTILLE{4.75}{-1.5}{-4}{-3}
   \POINTILLE{1}{-4.8}{-4.5}{2}
   \POINTILLE{-2.7}{-4.2}{-2}{4.5}
   \POINTILLE{-5}{0.3}{3.5}{3.5}
   \POINTILLE{-3.5}{3.5}{5}{0.3}
   \bSOM{0}{5}{2}{}
   \gSOM{2}{4.5}{}{}
   \gSOM{3.5}{3.5}{}{}
   \bSOM{4.5}{2}{}{}  \node[below] at (5.2,2.8) {1 };
   \gSOM{5}{0.3}{}{}
   \bSOM{4.75}{-1.5}{}{} \node[below] at (5.5,-0.8) {1};
   \gSOM{4}{-3}{}{}
   \gSOM{2.7}{-4.2}{}{}
   \bSOM{1}{-4.8}{}{1}
  \gSOM{-1}{-4.8}{}{}
  \bSOM{-2.7}{-4.2}{}{1}
  \gSOM{-4}{-3}{}{}
   \gSOM{-4.75}{-1.5}{}{}
   \bSOM{-5}{0.3}{}{} \node[below] at (-5.8,0.8) {1};
    \gSOM{-4.5}{2}{}{}
  \bSOM{-3.5}{3.5}{}{} \node[below] at (-4.1,4.5) {1};
     \gSOM{-2}{4.5}{}{}
  \node[rotate=-58] at (3.74,-2.5) {\textbf{{\large $>$}}};
  \node[rotate=170] at (-4.2,-1.58) {\textbf{{\large $>$}}};
\tobedone
\end{tikzpicture}

\caption{\label{fig:mapping 1}Construction of the mapping $g$ in the proof of Theorem~\ref{th:C(3a-1;1,a)} with $n=17$ and $a=6$.}
\end{center}
\end{figure}

Since $a$ is even, $g$ is  an independent broadcast on $\overrightarrow{C}(n;1,a)$ with cost  $$\sigma(g) = 2 + \left( \frac{a-4}{2}+1 \right) + \left( \frac{a-4}{2}+1 \right) +  \left( \frac{a-4}{2}+1 \right) = \frac{3a}{2} - 1 = \frac{n-1}{2}.$$
Hence $\beta_b(\overrightarrow{C}(n;1,a)) \geq \left\lfloor\frac{n}{2}\right\rfloor$.

By Observation~\ref{isomorphisme}, we then get
$\beta_b(\overrightarrow{C}(n;1,2)) = \beta_b(\overrightarrow{C}(n;1,-(2a-1))) = \left\lfloor\frac{n}{2}\right\rfloor.$
\end{enumerate}
This completes the proof.
\end{proof}

\begin{theorem}\label{th:C(n;1,2)}  
For every integer $n\ge 4$,
$$\beta_b(\overrightarrow{C}(n;1,2)) = \beta_b(\overrightarrow{C}(n;1,-(n-2)))  =  \diam (\overrightarrow{C}(n;1,2)) = \left\lfloor\frac{n}{2}\right\rfloor.$$
\end{theorem}

\begin{proof}
We first determine the diameter of the graph $\overrightarrow{C}(n;1,2)$. \\
For each vertex $v_i\in V(\overrightarrow{C}(n;1,2))$, $d(v_0,v_i)=\left\lceil\frac{i}{2}\right\rceil$. Then
 $$\diam (\overrightarrow{C}(n;1,2)) = e(v_0)= \max_{v_{i} \in V} d(v_{0},v_{i})=  \left\lceil\frac{n-1}{2}\right\rceil = \left\lfloor\frac{n}{2}\right\rfloor.$$
Therefore,  by Proposition~\ref{prop:lower bound 1}, $\beta_b(\overrightarrow{C}(n;1,2)) \geq \left\lfloor\frac{n}{2}\right\rfloor$. We now prove the reverse inequality. Let $f$ be an independent $\beta_b$-broadcast on $\overrightarrow{C}(n;1,2)$.  We consider two cases, depending on the value of $|V_f^+|$.
\begin{enumerate}
\item  $|V_f^+| = 1$. Then $V_f^+=\{v_i\}$ for some vertex $v_i$,
and thus
$$\beta_b(\overrightarrow{C}(n;1,2))=\sigma(f) = f(v_i) \le e(v_0) = \diam(\overrightarrow{C}(n;1,2)).$$

\item $|V_f^+| \ge 2$. In that case, each vertex $v\in V_f^+$ $f$-dominates $2f(v)+1$ vertices.
Moreover, each vertex is $f$-dominated at most once. This gives
$$2f(V_f^+) + |V_f^+| \leq n $$
and thus
$$\beta_b(\overrightarrow{C}(n;1,2)) = \sigma(f) = \sum_{v\in V_f^+}f(v) = f(V_f^+) \leq  \left \lfloor \frac{n - |V_f^+|}{2} \right \rfloor  \leq \left \lfloor \frac{n - 2}{2} \right \rfloor < \diam(\overrightarrow{C}(n;1,2)),$$
which is impossible. Hence $|V_f^+| = 1$ and  $$\beta_b(\overrightarrow{C}(n;1,2)) = \diam(\overrightarrow{C}(n;1,2)) = \left\lfloor\frac{n}{2}\right\rfloor.$$
By Observation~\ref{isomorphisme}, we then get
$\beta_b(\overrightarrow{C}(n;1,2)) = \beta_b(\overrightarrow{C}(n;1,2-n)) =  \diam (\overrightarrow{C}(2a;1,a)) = \left\lfloor\frac{n}{2}\right\rfloor.$
\end{enumerate}
\end{proof}
Since the circulant graphs $\overrightarrow{C}(2a-1;1,a)$ and $\overrightarrow{C}(2a-1;1,2)$ are isomorphic for every integer $a$, $a\ge 2$, (Figure~\ref{fig:iso} illustrates the two isomorphic graphs  $\overrightarrow{C}(9;1,2)$ and $\overrightarrow{C}(9;1,5)$) Theorem~\ref{th:C(n;1,2)} admits the following corollary.

\begin{figure}[h]
\begin{center}
\begin{tikzpicture}[scale=0.7]

 \foreach \i in {0,2,...,6} {
        \draw  (\i,0) -- (\i,2);}
   \foreach \i in {0,2} {
        \draw  (0,\i) -- (6,\i); }
             \draw  (6,2) -- (8,2);
               \foreach \x in {0,2,4,6,8}{
                    \node  at (\x ,2.5){ $v_{\x}$ } ;}
                      \foreach \x in {1,3,5,7}{
                        \node  at (\x-1 ,-0.5){ $v_{\x}$ } ;}

\foreach \i in {2,4,6,8} {
         \node[rotate=0] at (\i-0.5,2) {\textbf{{ $>$}}};}
 \foreach \i in {2,4,6} {
         \node[rotate=0] at (\i-0.5,0) {\textbf{{ $>$}}};}
    \foreach \i in {2,4,6,8} {
         \node[rotate=45] at (\i-0.3,1.7) {\textbf{{ $>$}}};}
\foreach \i in {2,4,6,8} {
         \node[rotate=-90] at (\i-2,0.5) {\textbf{{ $>$}}};}

      \node[rotate=210] at (0.5,2.3) {\textbf{{ $>$}}};
        \node[rotate=210] at (0.5,2.3) {\textbf{{ $>$}}};
         \node[rotate=140] at (0.4,-0.3) {\textbf{{ $>$}}};
          \node[rotate=160] at (0.7,1.76) {\textbf{{ $>$}}};
\node (A) at (8, 2) { };
\node  (B) at (0, 2) { };
\node  (C) at (6, 0) { };
\node  (D) at (0, 0) { };

\draw (A) to[bend right=30] (B);
\draw (D) to[out=-40, in=-70] (A);

\POINTILLE{0}{2}{6}{0}

 \foreach \i in {0,2,...,6} {
        \draw  (\i,0) -- (\i+2,2);
 }

 \foreach \x in {0,2,...,8}{
 \node [scale=0.7, draw, circle, fill=white] at (\x , 2){  } ;
}

 \foreach \x in {0,2,...,6}{
 \node [scale=0.7, draw, circle, fill=white] at (\x , 0){  } ;
}

 \foreach \i in {0,2,...,6} {
        \draw  (\i+11,0) -- (\i+11,2);}
   \foreach \i in {0,2} {
        \draw  (11,\i) -- (17,\i); }
             \draw  (17,2) -- (19,2);
               \foreach \x in {0,1,2,3,4}{
                    \node  at (\x + \x +11 ,2.5){ $v_{\x}$ } ;}
                      \foreach \x in {5,6,7,8}{
                        \node  at (\x + \x +1 ,-0.5){ $v_{\x}$ } ;}

\foreach \i in {2,4,6,8} {
         \node[rotate=0] at (\i+10,2) {\textbf{{ $>$}}};}
 \foreach \i in {2,4,6} {
         \node[rotate=0] at (\i+10.5,0) {\textbf{{ $>$}}};}
    \foreach \i in {2,4,6,8} {
         \node[rotate=45] at (\i+10.7,1.7) {\textbf{{ $>$}}};}
\foreach \i in {2,4,6,8} {
         \node[rotate=-90] at (\i+9,0.5) {\textbf{{ $>$}}};}

      \node[rotate=210] at (11.5,2.3) {\textbf{{ $>$}}};
        \node[rotate=210] at (11.5,2.3) {\textbf{{ $>$}}};
         \node[rotate=140] at (11.4,-0.3) {\textbf{{ $>$}}};
          \node[rotate=160] at (11.7,1.76) {\textbf{{ $>$}}};
\node (A) at (19, 2) { };
\node  (B) at (11, 2) { };
\node  (C) at (17, 0) { };
\node  (D) at (11, 0) { };

\draw (A) to[bend right=30] (B);
\draw (D) to[out=-40, in=-70] (A);

\POINTILLE{11}{2}{17}{0}

 \foreach \i in {0,2,...,6} {
        \draw  (\i+11,0) -- (\i+13,2);
 }

 \foreach \x in {0,2,...,8}{
 \node [scale=0.7, draw, circle, fill=white] at (\x+11 , 2){  } ;
}

 \foreach \x in {0,2,...,6}{
 \node [scale=0.7, draw, circle, fill=white] at (\x+11 , 0){  } ;
}
+-\node at (4 , -2){The directed graph $\overrightarrow{C}(9;1,2)$} ;
\node at (15 , -2){The directed graph $\overrightarrow{C}(9;1,5)$} ;

\end{tikzpicture}
\caption{\label{fig:iso} The isomorphic graphs  $\overrightarrow{C}(9;1,2)$ and $\overrightarrow{C}(9;1,5)$.}
\end{center}
\end{figure}

\begin{corollary}\label{cor:(C2a-1;1,a)}
For every integer $a\ge 2$,
$$\beta_b(\overrightarrow{C}(2a-1;1,a)) = \beta_b(\overrightarrow{C}(2a-1;1,-(a-1)))  =  \diam (\overrightarrow{C}(2a-1;1,a)) = \left\lfloor\frac{2a-1}{2}\right\rfloor.$$
\end{corollary}

\begin{theorem}\label{th:C(n;1,3)}
For every integer $n\ge 4$,
$$\beta_b(\overrightarrow{C}(n;1,3)) = \beta_b(\overrightarrow{C}(n;1,-(n-3)))= \left\lfloor\frac{n}{2}\right\rfloor.$$
\end{theorem}

\begin{proof}

Let $f$ be an independent $\beta_b$-broadcast on $\overrightarrow{C}(n;1,3)$. Each $f$-broadcast vertex $v_i\in V_f^+$, $f$-dominates all consecutive vertices from $v_i$ to $v_{i + 3f(v_i) - 2}$, ($3f(v_i)-1$ vertices). Therefore

$$ 3f(V_f^+) - \left|V_f^+\right| \leq n,$$
which gives

$$f(V_f^+) \leq \left\lfloor   \frac{n + \left|V_f^+\right|}{3} \right \rfloor.$$

Since $\left|V_f^{+}\right| \leq \dfrac{n}{2}$, we get

$$\sigma(f) = f(V_f^{+}) \leq \left\lfloor   \frac{n }{2} \right \rfloor.$$

We now prove the reverse inequality. For this, we construct a mapping $g$ from $V(\overrightarrow{C}(n;1,3))$ to $\{0,1,2\}$, considering two cases, depending on the parity of
$n$.

\begin{enumerate}
\item If $n$ is even, then we let $g(v_i)=1$ if
and only if $i$ is even.
Since $n$ is even, $g$ is an independent broadcast on $\overrightarrow{C}(n;1,3)$ of cost $ \sigma(g) = \frac{n}{2}$. Hence,
 $\beta_b(\overrightarrow{C}(n;1,3)) \geq \frac{n}{2}$.
\item If $n$ is odd, we let
$$ g(v_i) = \left\{
  \begin{array}{ll}
    1, & \mbox{if $i$ is even and  $0\leq i \leq n-7$,} \\ [1ex]
    2, & \mbox{if $i = n -5 $,} \\ [1ex]
    0, & \mbox{otherwise,}
  \end{array}
\right.$$
(see Figure~\ref{fig:mapping a = 3} for the case $n = 17$).

\begin{figure}[h]
\begin{center}
\begin{tikzpicture}[scale=0.5]

\draw(0,0)circle(5) ;
   \draw (0,6.5) arc (90: 0 : 6.5cm);
   \node[rotate=-90] at (6.5,0) {\textbf{{\large $>$}}};
 \LIGNE{0}{5}{4.5}{2}
 \LIGNE{4.5}{2}{4}{-3}
 \LIGNE{4}{-3}{-1}{-4.8}
\LIGNE{-1}{-4.8}{-4.75}{-1.5}
\LIGNE{-4.75}{-1.5}{-3.5}{3.5}
\POINTILLE{-3.5}{3.5}{2}{4.5}
   \bSOM{0}{5}{1}{}
   \gSOM{2}{4.5}{}{}
   \bSOM{3.5}{3.5}{}{} \node[below] at (4.5,4.1) {1 };
   \gSOM{4.5}{2}{}{}
   \bSOM{5}{0.3}{}{} \node[below] at (5.8,0.8) {1};
   \gSOM{4.75}{-1.5}{}{}
   \bSOM{4}{-3}{}{} \node[below] at (5,-2.6) {1};
   \gSOM{2.7}{-4.2}{}{}
   \bSOM{1}{-4.8}{}{1}
  \gSOM{-1}{-4.8}{}{}
  \bSOM{-2.7}{-4.2}{}{1}
  \gSOM{-4}{-3}{}{}
   \bSOM{-4.75}{-1.5}{}{} \node[below] at (-5.5,-0.8) {2};
   \gSOM{-5}{0.3}{}{}
    \gSOM{-4.5}{2}{}{}
  \gSOM{-3.5}{3.5}{}{}
     \gSOM{-2}{4.5}{}{}
 \node[rotate=-28] at (3.8,2.45) {\textbf{{\large $>$}}};
 \node[rotate=-95] at (4.1,-2) {\textbf{{\large $>$}}};
 \node[rotate=200] at (-0,-4.4) {\textbf{{\large $>$}}};
 \node[rotate=140] at (-4.05,-2.1) {\textbf{{\large $>$}}};
 \node[rotate=80] at (-3.7,2.5) {\textbf{{\large $>$}}};
\tobedone
\end{tikzpicture}
\caption{\label{fig:mapping a = 3}Construction of the mapping $f$ in the proof of Theorem~\ref{th:C(n;1,3)} with $n=17$.}
\end{center}
\end{figure}
Clearly $g$ is an independent broadcast on $\overrightarrow{C}(n;1,3)$ of cost $$\sigma(g) = 2 + \left( \frac{n-7}{2}+1 \right) = \frac{n-1}{2}.$$ Hence $\beta_b(\overrightarrow{C}(n;1,3)) \geq \left\lfloor\frac{n}{2}\right\rfloor$.\\

By Observation~\ref{isomorphisme}, we then get
$\beta_b(\overrightarrow{C}(n;1,3)) = \beta_b(\overrightarrow{C}(n;1,3-n)) =   \left\lfloor\frac{n}{2}\right\rfloor.$
\end{enumerate}
This completes the proof.
\end{proof}

\section{Independence broadcast number on $ \overrightarrow{C}(n;1,a)$, with $a \ge 4$ and $n>2a$  }\label{sec:Bb q>2}

In this section, we consider the oriented circulant graph $\overrightarrow{C}(n;1,a)$ with $n = qa + r$, $a \ge 4$, $q \ge 2$, and $0 \le r \le a-1$.
 The main objective of this section is to identify the maximum value that can be assigned to the broadcast vertices. This value allows us to give some bounds and some exact values of $\beta_b(\overrightarrow{C}(n;1,a))$.

\subsection{$\ell$-bounded optimal independent broadcasts}\label{sec:l-bounded}
Let $\overrightarrow{G}=\overrightarrow{C}(n;1,a)$ be an oriented circulant graph. Recall that we denote the vertices of $\overrightarrow{G}$ as $v_0, v_1, \dots, v_{n-1}$, and the subscripts are always considered modulo $n$. For a given broadcast $f$ in $\overrightarrow{G}$, we say that $f$ is \emph{$\ell$-bounded}, for some integer $\ell\ge 1$, if $f(v)\le \ell$ for every vertex $v$.\\
Our goal here is to identify some classes of oriented circulant graphs for which there exists an $\ell$-bounded independent $\beta_b$-broadcast. It is forth pointing that any bound $\ell$ cannot exceed the diameter of $\overrightarrow{G}$, since $f(v) \le \text{\diam}(\overrightarrow{G})$ for every vertex $v$. We now determine an upper bound of $\diam(\overrightarrow{G})$.

\begin{proposition}\label{prop:diam C(n;1,a)}
If $n$, $a$, $q$ and $r$ are four integers, with $n = qa + r$, $a \ge 4$, $q \ge 2$, and $0 \le r \le a-1$ then
$$ \diam (\overrightarrow{C}(n;1,a)) \leq q +a-2.$$
\end{proposition}

\begin{proof}
For each vertex $v_{i}\in V=V(\overrightarrow{C}(n;1,a))$, there exists a path connecting $v_0$ to $v_i$, defined by
\begin{center}
$P_i = \{v_{0}, v_{a}, v_{2a}, \ldots, v_{pa}, v_{pa+1},\ldots,v_{i}\}$, where $p=\left\lfloor \frac{i}{a}\right\rfloor $,
\end{center}
which gives

\begin{center}
$d(v_{0},v_{i}) \leq  p + i - pa$, where $p + i - pa$ is the length of the path $P_i$.
\end{center}

Since  $\diam(\overrightarrow{C}(n;1,a)) = e(v_0) $, we get

\begin{align*}
\diam (\overrightarrow{C}(n;1,a))
  & = \max_{v_i \in V} \left(d(v_{0},v_{i})  \right) \le \max_{0 \le i \le n-1} \left( \left\lfloor \frac{i}{a}\right\rfloor + i - \left\lfloor \frac{i}{a}\right\rfloor a   \right)
  \\[2ex]
  & \leq \max \{ q-1 + a-1 ; q + r - 1 \}  = q+a-2.
 \end{align*}
This completes the proof.
\end{proof}
In order to prove that the diameter bound in Proposition~\ref{prop:diam C(n;1,a)} is achieved, let us start by establishing a result relating to the distance between any two vertices of $\overrightarrow{C}(n;1,a)$, with $n=qa$, $a\ge 4$ and $ q \ge 3$ or $n=qa+r$ with $q\geq 2$, $1 \le r \le a-1$ and $ 4 \le a \leq r + q + 1.$
\begin{proposition}\label{prop:Distance1}
Let $n, a, q$ and $r$ be four integers such that $n=qa$, $a\ge 4$ and $ q \ge 3$ or $n=qa+r$ with $q\geq 2$, $1 \le r \le a-1$ and $ 4 \le a \leq r + q + 1.$ For every two vertices $v_{i},v_{j}$, $i<j$,  of $\overrightarrow{C}(n;1,a)$, we have
$$ d(v_{i},v_{j}) = \left\lfloor \frac{j-i}{a}\right\rfloor (1-a) + j -  i.$$
\end{proposition}
\begin{proof}
Let $\overrightarrow{C}(n;1,a)$ be an oriented circulant graph and $v_i$, $v_{j}$ two vertices of $\overrightarrow{C}(n;1,a)$ such that $i<j$. For any $v_{i}-v_{j}$ path $P$,  there exists an integer $\ell$ such that $$ P = P_\ell = \{ v_{i},v_{i+a}, \ldots,v_{ i+\ell a},v_{i+\ell a + 1 }, \ldots, v_{j} \}$$ and
$$ d(v_{i},v_{j}) = \min_{\ell \geq 0}  \vert E(P_\ell)  \vert = \min_{\ell \geq 0} \ell + (j- (i + \ell a)) $$
where $\vert E(P_\ell)\vert $ is the length of the path $P_\ell$.

It is easy to see that the case $r=0$, the shortest $v_i - v_j$  path is $$P = \{ v_{i},v_{i+a}, \ldots ,v_{ i+ \left\lfloor \frac{j-i}{a}\right\rfloor a},v_{i+\left\lfloor \frac{j-i}{a}\right\rfloor a + 1}, \ldots , v_j \}.$$
 Assume now $r \neq 0$. We consider four cases, depending on the value of $\ell$ to prove that $$ d(v_{i},v_{j}) = \min_{\ell \geq 0}  \vert E(P_\ell)  \vert =  \left\lfloor \frac{j-i}{a}\right\rfloor (1-a) + j -  i.$$
\begin{enumerate}
\item If $ 0 \le \ell < \left\lfloor \frac{j-i}{a}\right\rfloor $, then
 $$ \vert E(P_\ell)  \vert  = \ell + j - (i+ \ell a ) = \ell(1 - a) + j - i.$$
 Since $ 2 < a $, we get
$$ \vert E(P_\ell)  \vert  >  \left\lfloor \frac{j-i}{a}\right\rfloor(1 - a) + j - i.$$

\item If $\ell = \left\lfloor \frac{j-i}{a}\right\rfloor  $, then

$P_\ell = \{v_{i},v_{i+a},v_{i+2a}, \ldots ,v_{i+\left\lfloor \frac{j-i}{a}\right\rfloor a}, v_{i+\left\lfloor \frac{j-i}{a}\right\rfloor a+1}, \ldots , v_{j} \}$. Therefore,
$$ \vert E(P_\ell)  \vert =  \left\lfloor \frac{j-i}{a}\right\rfloor (1-a) + j -  i. $$

\item If $ \left\lfloor \frac{j-i}{a}\right\rfloor + 1 \le \ell \le q $, then $ j < \ell a $ and

$$P_\ell = \{v_{i},v_{i+a},v_{i+2a}, \ldots , v_{\ell a}, v_{\ell a +1 }, \ldots, v_{i-1}, v_{i}, v_{i+1}, \ldots, v_{j}\}.$$

It is clearly seen that the vertex $v_{i}$ appears twice in the path $P_\ell$, which implies that this path is not minimal.

\item If $ q+1 \le \ell $, then
 $$P_\ell = \{v_{i},v_{i+a},v_{i+2a}, \ldots , v_{i+qa}, v_{i+(q+1)a}, \ldots , v_{i+\ell a},v_{i+\ell a + 1}, v_{i+\ell a + 2}, \ldots , v_{j}\}$$

Let  $P'_\ell$ and $P'_{\ell^{'}}$ be two paths defined by
\begin{center}
$ P'_\ell= \{v_{i},v_{i+a},v_{i+2a}, \ldots , v_{i+(q+1)a} = v_{i+a-r}  \} $ and $ P'_{\ell^{'}}= \{v_{i},v_{i+1},v_{i+2}, \ldots , v_{i+a-r} \}$

\end{center}

Since $a \le q+r+1$, we get  $ \vert E(P'_\ell)  \vert   \geq  \vert E(P'_{\ell^{'}})  \vert $. Moreover,  $P'_\ell$ is a sub-path of $P_\ell$,
then there exists a path $P_{\ell^{'}}$  such that $\ell^{'} \le q$ and  $ \vert E(P_\ell)  \vert \geq \vert E(P_{\ell^{'}})\vert$.
It follows
$$ \vert E(P_\ell)  \vert \geq \vert E(P_{\ell^{'}})\vert \geq \left\lfloor \frac{j-i}{a}\right\rfloor(1 - a) + j - i $$
\end{enumerate}

In all cases, we get $\vert E(P_\ell)  \vert \geq \left\lfloor \frac{j-i}{a}\right\rfloor(1 - a) + j - i $, and the equality is reached in the case $\ell = \left\lfloor \frac{j-i}{a}\right\rfloor$, which implies that
$$ d(v_i,v_{j}) =  \min_{\ell \geq 0}  \vert E(P_\ell)  \vert =  \left\lfloor \frac{j-i}{a}\right\rfloor(1 - a) + j - i.$$ This completes the proof.
\end{proof}

\begin{proposition}\label{prop:Diam}
Let $n, a, q$ and $r$ be four integers such that $n=qa$ with $a\ge 4$ and $ q \ge 3$ or $n=qa+r$ with $q\geq 2$, $1 \le r \le a-1$ and $ 4 \le a \leq r + q + 1.$ For every two vertices $v_{i},v_{j}$, $i<j$, of $\overrightarrow{C}(n;1,a)$
we have
$$ \diam (\overrightarrow{C}(n;1,a)) = q + a - 2.$$
\end{proposition}

\begin{proof}
Let $V=V(\overrightarrow{C}(n;1,a))$. From {Proposition}~\ref{prop:Distance1}, and since $\diam (\overrightarrow{C}(n;1,a)) = e(v_0)$, we get

$$\diam (\overrightarrow{C}(n;1,a)) = \max_{v_{i} \in  V } d(v_{0},v_{i})=  \max_{1 \le i \le n-1 } \left(\left\lfloor \frac{i}{a}\right\rfloor(1 - a) + i\right) = q  + a - 2.$$

This value is given by the distance between $v_0$ and $v_{qa-1}$. This completes the proof.
\end{proof}

Now, we will begin the process of identifying two various $\ell$-bounded independent broadcasts on circulant graphs of the form $\overrightarrow{C}(qa + r; 1, a)$, where $a\ge 4$, $q \geq 2$  and $0 \le r \le a-1$. We will first consider the case where $a \leq q+r+1$.

\begin{lemma}\label{lem:broadcast-at-most-a-1-a<q+r}
Let $n, a, q$ and $r$ be four integers such that $n=qa+r$. If $4 \le a \leq q+r+1$ and $0 \le r \le a-2$, or $a \le q + 1 $ and $r=a-1$, then  $\overrightarrow{C}(n;1,a)$ admits an $(a-1)$-bounded $\beta_b$-broadcast.

\end{lemma}

\begin{proof}
It is enough to prove that for every independent broadcast $f$ on $\overrightarrow{C}(n;1,a)$,
there exists an independent broadcast $g$ on $\overrightarrow{C}(n;1,a)$ such that $\sigma(g)\ge\sigma(f)$
and $g(v)\leq a-1$ for every vertex $v \in V_{g}^+$.\\
Let $f$ be any independent broadcast on $\overrightarrow{C}(n;1,a)$ and $g$ be the mapping from $V(\overrightarrow{C}(n;1,a))$ to $\{0,1 \ldots, a-1\}$ defined as follows
(the construction of the mapping $g$ is illustrated in Figure~\ref{fig:mapping-g},
not all $a$-edges being drawn).

\begin{enumerate}
\item  If $v_i$ is an $f$-broadcast vertex such that $a \le f(v_{i}) \le 2a-3 $, then we let
$$g(v_j) =
\left\{
   \begin{array}{ll}
   a-2 & \text{if } j=i, \\ [1ex]
   d & \text{if $j=i+(a-1)$  }, \\ [1ex]
   1 & \text{if $j=i+2a-2$  }, \\ [1ex]
   \end{array}
\right.
$$
where $d = f(v_{i})- (a - 1)$ (see Figure~\ref{fig:mapping-g}(a)).
\item  If $v_i$ is an $f$-broadcast vertex such that $p(a-1)\le f(v_{i}) \le (p+1)(a-1)-1 $, for some $p \geq 2$,  then we let
$$g(v_j) =
\left\{
   \begin{array}{ll}
   a-2 & \text{if } i\le j \le i + (4p-5)(a-1), \\ [1ex]
   & \hskip 0.5cm \text{and} (j-i)\mod (a-1) \text{ is odd},
   \end{array}
\right.
$$
(see Figure~\ref{fig:mapping-g}(b)).
\begin{figure}[h]
\begin{center}
\begin{tikzpicture}[scale=0.5]
   \node[above] at (0,0.2){$f$};
   \node[below] at (0,-0.2){$g$};

   \POINTILLE{0}{0}{1}{0}   \LIGNE{1}{0}{27}{0}   \POINTILLE{27}{0}{28}{0}
   \foreach \k in {2,8,14,20} \draw[thick] (\k,0) .. controls (\k+1,2) and (\k+5,2) .. (\k+6,0);
   \foreach \k in {1,2,...,27} \gSOM{\k}{0}{}{};
   \SOM{1}{0}{}{ }
   \bSOM{2}{0}{}{4}   \node[above] at (2,0.5) {7};
   \bSOM{7}{0}{}{2}
   \bSOM{12}{0}{}{1}
   \SOM{17}{0}{}{ }
   \SOM{18}{0}{}{ }
   \SOM{19}{0}{}{ }
   \SOM{23}{0}{}{ }
   \SOM{24}{0}{}{ }
    \SOM{25}{0}{}{ }
    \SOM{27}{0}{}{ }
   \node[below] at (14,-2){(a) Item 1: $f(v_i)=7$ and $a=6$ (so that $d=2$)};
\end{tikzpicture}

\vskip 0.5cm
\begin{tikzpicture}[scale=0.5]
   \node[above] at (0,0.2){$f$};
   \node[below] at (0,-0.2){$g$};
   \POINTILLE{0}{0}{1}{0}   \LIGNE{1}{0}{27}{0}   \POINTILLE{27}{0}{28}{0}
   \foreach \k in {2,6,10,14,18,22} \draw[thick] (\k,0) .. controls (\k+1,1.8) and (\k+3,1.8) .. (\k+4,0);
   \foreach \k in {5,8,9,11,12,15} \draw[dotted] (\k,0) .. controls (\k+1,1.8) and (\k+3,1.8) .. (\k+4,0);
   \foreach \k in {1,2,...,19} \gSOM{\k}{0}{}{};
 \foreach \k in {20,21,...,27} \SOM{\k}{0}{}{};
   \bSOM{2}{0}{}{2} \node[above] at (2,0.5) {7};
   \bSOM{5}{0}{}{2}
   \bSOM{8}{0}{}{2}
   \bSOM{11}{0}{}{2}
    \SOM{14}{0}{}{ }
    \SOM{17}{0}{}{ }
     \SOM{18}{0}{}{ }
   \node[below] at (14,-2){(b) Item 2: $f(v_i)=7$ and $a=4$ (so that $p=2$)};
\end{tikzpicture}

\caption{\label{fig:mapping-g}Construction of the mapping $g$ in the proof of Lemma~\ref{lem:broadcast-at-most-a-1-a<q+r}.}
\end{center}
\end{figure}
\item For every other vertex $v_k$, we let $g(v_k)=f(v_k)$
\end{enumerate}

We now prove that in  all cases, $g$ is an independent broadcast on $\overrightarrow{C}(n; 1, a)$. For that, we first prove the following claim.

\begin{claim}\label{cl:distance-2a+r}
For every vertex $v_j$ whose $g$-value is set to~1, $d$ or $a-2$ in Item 1, or is set to $a-2$ in Item 2, we have
$d(v_i,v_j)\le f(v_i)-g(v_j)$.
\end{claim}
\begin{proof}

In Item~1, thanks to {Proposition}~\ref{prop:Distance1}, we have
$$d(v_i,v_{i+a-1}) = a-1 \le f(v_i) - d = f(v_i)-g(v_{i+a-1}),$$
and
$$d(v_i,v_{i+2a-2}) = a-1 \le f(v_i) - 1 = f(v_i)-g(v_{i+2a-2}).$$

According to {Proposition}~\ref{prop:Distance1}, the distance between vertex $v_i$ and all vertices whose $g$-value could be set to $a-2$ in Item~2 is the same, and it is given by $d(v_i,v_j)= a-1 $. Moreover, since $p \geq 2$, we get
 $$d(v_i,v_j) \le  (p-1)(a-1) < p(a-1) - (a-2) = f(v_i) - g(v_j).$$
This concludes the proof of the claim.
\end{proof}

Thanks to Claim~\ref{cl:distance-2a+r}, and since  $f $ is an independent broadcast on $\overrightarrow{C}(n; 1, a)$, no $g$-broadcast vertex $v_i$ in some item dominates another $g$-broadcast vertex $v_j$  in another item. Therefore, to prove that $g$ is indeed an independent  broadcast on $\overrightarrow{C}(n; 1, a)$, it remains to prove the following claim.
\begin{claim}\label{cl:distance-general}
For every two $g$-broadcast vertices $v_i$ and $v_j$ in some item, we have
$$d(v_{i},v_{j}) >  g(v_{i})$$
\end{claim}
\begin{proof}

In Item~1, and  from {Proposition}~\ref{prop:Distance1}, we have
\[\begin{array}{l}
   d(v_i,v_{i+a-1}) = d(v_i,v_{i+2a-2}) = a-1 > g(v_i), \\ [1ex]
   d(v_{i+a-1},v_{i+2a-2}) = a-1 > d = g(v_{i+a-1}), \\ [1ex]
   d(v_{i+a-1},v_{i}) = \min \{q + r, q \}= q > d = g(v_{i+a-1}), \mbox{ and}  \\ [1ex]
   \min \{d(v_{i+2a-2},v_{i}), d(v_{i+2a-2},v_{i+a-1})\} > 1 = g(v_{i+2a-2}).
   \end{array}
\]

Now, in Item~2, and from {Proposition}~\ref{prop:Distance1}, we have for every two vertices $v_{j}$, $v_{j'}$ whose $g$-value is set to~$a-2$
 $$d(v_{j},v_{j'}) = a-1 > g(v_j).$$

This concludes the proof of the claim.
\end{proof}

To finish the proof, we only need to show that we have $\sigma(g) \geq \sigma(f)$. For every $v_i$ in Item 1, we have $g(v_i)+g(v_{i+a-1})+g(v_{i+2a-2}) = a-1+d = f(v_i)$. In Item 2, the number of vertices set to $a-2$ is $ n_1 = 4(p-1)$. We consider three cases, depending on the values of $a$ and $p$.

\begin{enumerate}
\item If $p=2$ and $a\geq 4$, then  $ n_1(a-2) = 4(a-2) \geq 3(a-1)-1 \geq f(v_i)$.

\item If $a=2$ and $p\geq 2$, then  $ n_1(a-2) = 8(p-1) \geq 3(p+1)-1 \geq f(v_i)$.

\item If $p \geq 3$ and $a \geq 5$, then  $ n_1(a-2) = 4(p-1)(a-2)= (2p-2)(2a-4) \geq (p+1)(a-1) > f(v_i)$.
\end{enumerate}

In all cases, we have $n_1(a-2) \geq f(v_i)$. We thus have $\sigma(g) \geq \sigma(f)$, as required. This completes the proof.
\end{proof}

Consider the case of circulant graphs $\overrightarrow{C}(n;1,a)$ when $a $ divides $n$.

\begin{lemma}\label{lem:broadcast-at-most-qa}
If $n$, $q$ and $a$  are three integers such that $ n= qa $, $ 3 \le q $,
then $\overrightarrow{C}(n;1,a)$ admits an $(a-1)$-bounded $\beta_b$-broadcast if $a - 1 \le q $, and an $q$-bounded $\beta_b$-broadcast if $a - 1 > q.$
\end{lemma}

\begin{proof}
Before starting the proof, let us note that the oriented circulant graph $\overrightarrow{C}(qa;1,a)$ can be illustrated as a modified grid with additional arcs or again, as an infinite grid created by connected repetitions of the grid $P_a \square P_q$ (Figure~\ref{fig:Ci} illustrates the circulant graph $\overrightarrow{C}(qa;1,a)$, with $q=3$ and $a=7$). For an independent broadcast $f$ on $\overrightarrow{C}(qa;1,a)$ this illustration allows us to have a better view of the $f$-dominated vertices of the graph.

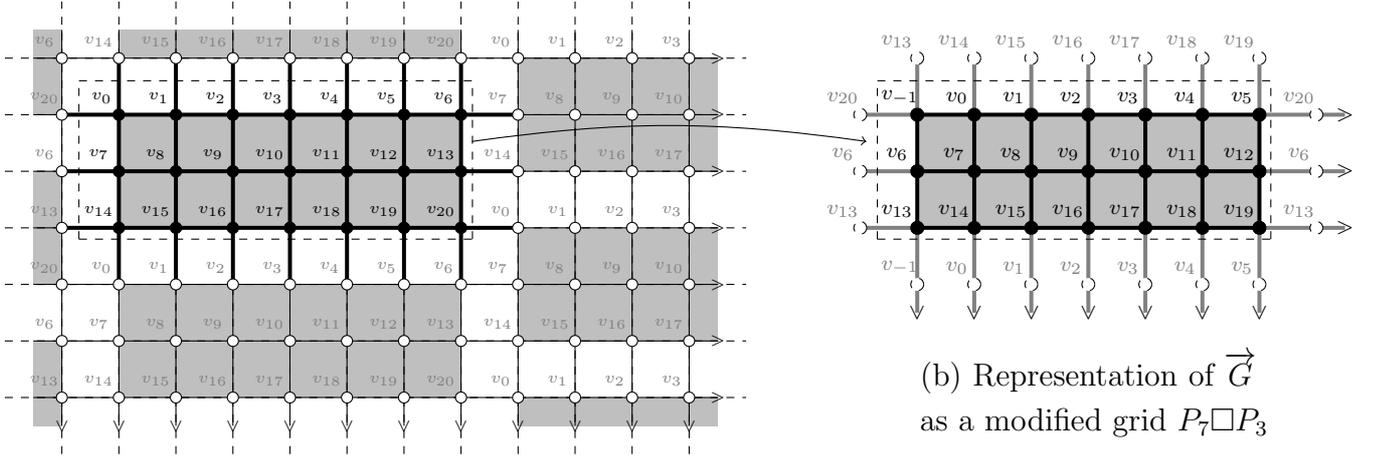
\begin{figure}[h]
\begin{center}
\begin{tikzpicture}[scale=0.75]

\begin{tiny}
\fill[lightgray] (2,3) rectangle (8,5);
\fill[lightgray] (2,6.5) rectangle (8,6);
\fill[lightgray] (1,2) rectangle (0.5,4);
\fill[lightgray] (1,-0.5) rectangle (0.5,1);
\fill[lightgray] (1,6.5) rectangle (0.5,5);
\fill[lightgray] (2,0) rectangle (8,2);
\fill[lightgray] (12.5,4) rectangle (9,6);
\fill[lightgray] (12.5,1) rectangle (9,3);
\fill[lightgray] (12.5,1) rectangle (9,3);
\fill[lightgray] (12.5,-0.5) rectangle (9,0);
 \foreach \j in {0,1,...,6} {
        \draw [dashed] (0,\j) -- (13,\j);
    }
  \foreach \i in {1,2,...,12} {
        \draw [dashed] (\i,-1) -- (\i,7);
        }
 \foreach \j in {0,1,...,6} {
        \draw  (0.5,\j) -- (12.5,\j);
    }
  \foreach \i in {1,2,...,12} {
        \draw  (\i,-0.5) -- (\i,6.5);

}
 \foreach \j in {3,4,5} {
        \draw[line width=1.5pt]  (1,\j) -- (9,\j);
    }
  \foreach \i in {2,3,...,8} {
       \draw[line width=1.5pt]  (\i,2) -- (\i,6);
}
    \foreach \x in {1,2,...,12} {
   \node [rotate=-90] at (\x ,-0.5){ $>$ } ;
}
   \foreach \x in {0,1,...,6} {
   \node  at ( 12.5,\x){ $>$ } ;
}
 \foreach \x in {1,2,...,12}
 \foreach \y in {0,1,...,6} {
        \node[scale=0.6, draw, circle, fill=white] (point\x) (point\x) at (\x,\y ) {};
    }
    \foreach \x in {2,3,...,8}
 \foreach \y in {3,4,5} {
        \node[scale=0.6, draw, circle, fill=black] (point\x) (point\x) at (\x,\y ) {};
    }

  \draw [dashed] (1.3,2.8) -- (8.2,2.8);
  \draw [dashed] (1.3,5.6) -- (8.2,5.6);
  \draw [dashed] (1.3,2.8) -- (1.3,5.6);
  \draw [dashed] (8.2,2.8) -- (8.2,5.6);

 \foreach \x in {0,1,...,6}{
 \node  at (\x + 1.7,5.3){ $v_{\x}$ } ;
 \node[gray] at (\x + 1.7,2.3){ $v_{\x}$ } ;
 }
 \foreach \x in {7,8,...,13}{

 \node at (\x -5.35,4.3){ $v_{\x}$ } ;
 \node [gray] at (\x -5.35,1.3){ $v_{\x}$ } ;
 }
  \foreach \x in {14,15,...,20}{
  \node [gray]   at (\x -12.35,6.3){ $v_{\x}$ } ;
 \node at (\x -12.35,3.3){ $v_{\x}$ } ;
 \node [gray] at (\x -12.35,0.3){ $v_{\x}$ } ;
 }

  \foreach \x in {0,1,...,3}{
 \node[gray]  at (\x + 8.7,6.3){ $v_{\x}$ } ;
 \node[gray] at (\x + 8.7,3.3){ $v_{\x}$ } ;
 \node[gray] at (\x + 8.7,0.3){ $v_{\x}$ } ;
 }
 \foreach \x in {7,8,...,10}{

 \node[gray] at (\x +1.65,5.3){ $v_{\x}$ } ;
 \node [gray] at (\x +1.65,2.3){ $v_{\x}$ } ;
 }
  \foreach \x in {14,15,...,17}{

 \node [gray]at (\x -5.35,4.3){ $v_{\x}$ } ;
 \node [gray] at (\x -5.35,1.3){ $v_{\x}$ } ;
 }

  \foreach \y in {0,3}{

 \node [gray]at (0.7 , \y + 2.3){ $v_{20}$ } ;
 \node [gray]at (0.7 , \y + 0.3){ $v_{13}$ } ;
 \node [gray]at (0.7 , \y + 1.3){ $v_{6}$ } ;
 }
  \node [gray]at (0.7 , 6.3){ $v_{6}$ } ;

\end{tiny}

\begin{scriptsize}

\fill[lightgray] (16,3) rectangle (22,5);

 \foreach \x in {-1,0,...,5}{

 \node  at (\x + 16.7,5.3){ $v_{\x}$ } ;
 \node [gray]at (\x + 16.7,2.3){ $v_{\x}$ } ;
 }
 \foreach \x in {6,...,12}{

 \node at (\x + 9.65,4.3){ $v_{\x}$ } ;
 }
  \foreach \x in {13,14,...,19}{

 \node at (\x + 2.65,3.3){ $v_{\x}$ } ;
  \node [gray]at (\x + 2.65,6.3){ $v_{\x}$ } ;
 }
  \foreach \y in {3}{
 \node [gray]at (14.7 , \y + 2.3){ $v_{20}$ } ;
 \node [gray]at (14.7 , \y + 0.3){ $v_{13}$ } ;
 \node [gray]at (14.7 , \y + 1.3){ $v_{6}$ } ;
 }
  \foreach \y in {3}{
 \node [gray]at (22.7 , \y + 2.3){ $v_{20}$ } ;
 \node [gray]at (22.7 , \y + 0.3){ $v_{13}$ } ;
 \node [gray]at (22.7 , \y + 1.3){ $v_{6}$ } ;
 }
  \foreach \x in {16,18,...,22}
 \foreach \y in {3,4,5} {
        \node[scale=0.6, draw, circle, fill=black] (point\x) (point\x) at (\x,\y ) {};
}
 \foreach \j in {3,4,5} {
        \draw[line width=1.5pt]  (16,\j) -- (22,\j);
    }
  \foreach \i in {16,...,22} {
       \draw[line width=1.5pt]  (\i,3) -- (\i,5);
}
 \foreach \j in {3,4,5} {
        \draw[gray,line width=1.5pt]  (15,\j) -- (16,\j);
         \draw[gray,line width=1.5pt]  (22,\j) -- (23.5,\j);
    }
\foreach \i in {16,...,22} {
       \draw[gray,line width=1.5pt]  (\i,5) -- (\i,6);
          \draw[gray,line width=1.5pt]  (\i,1.5) -- (\i,3);
}
     \foreach \x in {16,17,...,22}
 \foreach \y in {3,4,5} {
        \node[scale=0.6, draw, circle, fill=black] (point\x) (point\x) at (\x,\y ) {};
}
     \foreach \x in {16,17,...,22}
 \foreach \y in {2,6} {
        \node[scale=0.6, draw,dashed, circle, fill=white] (point\x) (point\x) at (\x,\y ) {};
}
 \foreach \x in {15,23}
 \foreach \y in {3,4,5} {
        \node[scale=0.6, draw,dashed, circle, fill=white] (point\x) (point\x) at (\x,\y ) {};
}
   \foreach \x in {16,17,...,22} {
   \node [rotate=-90] at (\x ,1.5){ $>$ } ;
}
   \foreach \x in {3,4,5} {
   \node  at ( 23.5,\x){ $>$ } ;
}
\end{scriptsize}
\node (A) at (8, 4.5) { };
  \node  (B) at (15.3, 4.5) { };
  \draw[->, bend left=8] (A) to (B);
    \draw [dashed] (15.3,2.8) -- (22.2,2.8);
  \draw [dashed] (15.3,5.6) -- (22.2,5.6);
  \draw [dashed] (15.3,2.8) -- (15.3,5.6);
  \draw [dashed] (22.2,2.8) -- (22.2,5.6);
 \node[below] at (19,1){ (b) Representation of $\overrightarrow{G}$} ;
 \node[below] at (19.1,0){  as a modified grid $P_7 \square P_3$ } ;
\end{tikzpicture}
\caption{\label{fig:Ci} The circulant graph $\overrightarrow{C}(qa;1,a)$ with $q=3$ and $a=7$}
\end{center}
\end{figure}

Let $f$ be an independent broadcast on $\overrightarrow{C}(qa;1,a)$. In order to prove the statement, we have to establish that there exists an independent broadcast $g$ on $C(qa; 1, a)$ satisfying $\sigma(g) \geq \sigma(f)$ and, for every vertex $v \in V_{g}^{+}$,  $g(v) \leq a-1$ if $a-1 \leq q$ or $g(v) \leq q$ if $q < a-1$. We consider two cases, according to the  maximum value between $a$ and $q-1$, to define the mapping $g$ from $V(\overrightarrow{C}(qa;1,a))$ to $\{0,1 \ldots, \min \{a-1, q\}\}$.
\begin{enumerate}
\item   $a - 1 \le q$. In this case, the result directly follows from Lemma~\ref{lem:broadcast-at-most-a-1-a<q+r} (The same construction for the function $g$ can be applied here as well).
\item $q < a-1$.
\begin{enumerate}
\item  If $v_i$ is an $f$-broadcast vertex such that $ q+1 \le f(v_{i}) \le 2q-1 $, then we let
$$g(v_j) =
\left\{
   \begin{array}{ll}
   q-1 & \text{if } j=i, \\ [1ex]
   d & \text{if $j= i + (q-k)a + k$ and  $1 \le k \le d + 1$}, \\ [1ex]
   \end{array}
\right.
$$
where $d = f(v_i)-q  $ (see Figure~\ref{fig:mapping-g-1}(a)).

\item $v_i$ is an $f$-broadcast vertex such that $pq\le f(v_{i}) \le (p+1)q-1 $ and $p \geq 2$,  then we let
\begin{center}
$g(v_j) = q-1$ if $j=i$ or ($d  \le j \le d + (p-2)q$ and $(j - d) \equiv 0 \mod (q)$)
\end{center}
where $d = i+(q-k)a+k $ and $ 1 \le k \le q $ (see Figure~\ref{fig:mapping-g-1}(b))

\item For every other vertex $v_i$, we let $g(v_i)=f(v_i)$
\end{enumerate}
\end{enumerate}
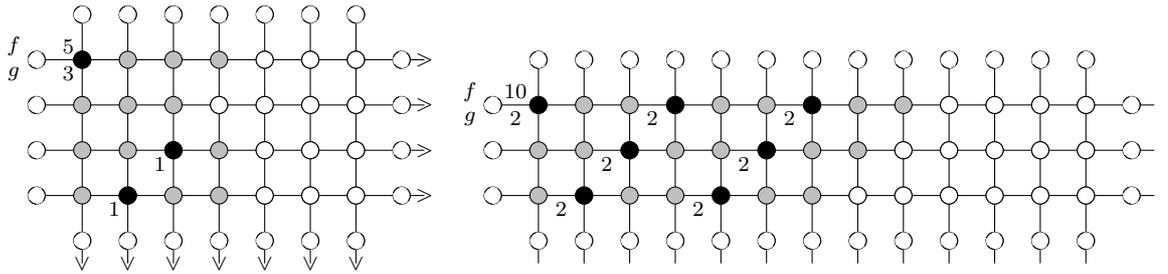
\begin{figure}[h]
\begin{center}
\begin{tikzpicture}[scale=0.6]
\begin{scriptsize}
 \foreach \i in {0,1,...,6} {
        \draw  (\i,-1.5) -- (\i,4);
 }
   \foreach \i in {0,1,2,3} {
        \draw  (-1,\i) -- (7.5,\i);
 }
 \foreach \x in {-1,...,7}{
 \node [scale=0.8, draw, circle, fill=lightgray] at (\x , 0){  } ;
 \node [scale=0.8, draw, circle, fill=lightgray] at (\x  , 1){  } ;
 \node [scale=0.8, draw, circle, fill=lightgray] at (\x  , 2){  } ;
  \node [scale=0.8, draw, circle, fill=lightgray] at (\x  , 3){  } ;
 }
 \foreach \x in {0,...,6}{
 \node [scale=0.8, draw, circle, fill=lightgray] at (\x , 4){  } ;
 \node [scale=0.8, draw, circle, fill=lightgray]  at (\x , -1){  } ;
 }
\foreach \x in {0,...,6}{
 \node [scale=0.8, draw,dashed, circle, fill=white]  at (\x , 4){  } ;
\node [scale=0.8, draw,dashed, circle, fill=white]  at (\x , -1){  } ;
}
\foreach \x in {0,1,2,3}{
 \node [scale=0.8, draw, circle,dashed, fill=white]  at (-1,\x ){  } ;
\node [scale=0.8, draw, circle,dashed, fill=white]  at (7,\x){  } ;
}
   \foreach \x in {0,1,...,6} {
   \node [rotate=-90] at (\x ,-1.5){ $>$ } ;
}
   \foreach \x in {0,1,2,3} {
   \node  at ( 7.5,\x){ $>$ } ;
}
    \node at (-1.5,3.3){$f$} ;
    \node at (-1.5,2.7){ $g$} ;
    \node [scale=0.8, draw, circle, fill=black] at (0 , 3){  } ;
     \node at (-0.3,3.3){$5$} ;
     \node at (-0.3,2.7){$3$} ;
    \node [scale=0.8, draw, circle, fill=black] at (1 , 0){  } ;
    \node at (1.7,0.7){$1$} ;
    \node [scale=0.8, draw, circle, fill=black] at (2 , 1){  } ;
    \node at (0.7,-0.3){$1$} ;

   \foreach \x in {4,5,6}{
    \node [scale=0.8, draw, circle, fill=white] at (\x , 3){  } ;
    }
    \foreach \x in {3,4,5,6}{
    \node [scale=0.8, draw, circle, fill=white] at (\x , 2){  } ;
    }
   \foreach \x in {4,5,6}{
    \node [scale=0.8, draw, circle, fill=white] at (\x ,1){  } ;
    }
  \foreach \x in {4,5,6}{
    \node [scale=0.8, draw, circle, fill=white] at (\x ,0){  } ;
    }
\end{scriptsize}
\node[below] at (3,-2){(a) : $f(v_i)= 5$, $q = 4$, $a=7$ (so that $d=1$)} ;
\begin{scriptsize}
\foreach \i in {10,...,22} {
        \draw  (\i,-1.5) -- (\i,3);
 }
   \foreach \i in {0,1,2} {
        \draw  (9,\i) -- (23.5,\i);
 }
 \foreach \x in {9,...,23}{
 \node [scale=0.8, draw, circle, fill=lightgray] at (\x , 0){  } ;
 \node [scale=0.8, draw, circle, fill=lightgray] at (\x  , 1){  } ;
 \node [scale=0.8, draw, circle, fill=lightgray] at (\x  , 2){  } ;
 }
 \foreach \x in {10,...,22}{
 \node [scale=0.8, draw, circle, fill=lightgray] at (\x , 3){  } ;
 \node [scale=0.8, draw, circle, fill=lightgray]  at (\x , -1){  } ;
 }
\foreach \x in {10,...,22}{
 \node [scale=0.8, draw, circle,dashed, fill=white]  at (\x , 3){  } ;
\node [scale=0.8, draw, circle,dashed, fill=white]  at (\x , -1){  } ;
}
\foreach \x in {0,1,2}{
 \node [scale=0.8, draw,dashed, circle, fill=white]  at (9,\x ){  } ;
\node [scale=0.8, draw,dashed, circle, fill=white]  at (23,\x){  } ;
}
    \node at (8.5,2.3){$f$} ;
    \node at (8.5,1.7){ $g$} ;
     \node at (9.5,2.3){$10$} ;
    \foreach \x in {10,13,16}{
     \node at (\x-0.5,1.7){$2$} ;
    \node [scale=0.8, draw, circle, fill=black] at (\x , 2){  } ;
    }
    \foreach \x in {12,15}{
     \node at (\x-0.5,0.7){$2$} ;
    \node [scale=0.8, draw, circle, fill=black] at (\x , 1){  } ;
    }
   \foreach \x in {11,14}{
     \node at (\x-0.5,-0.3){$2$} ;
    \node [scale=0.8, draw, circle, fill=black] at (\x ,0){  } ;
    }
    \foreach \x in {19,20,21,22}{
    \node [scale=0.8, draw, circle, fill=white] at (\x ,2){  } ;
    }
    \foreach \x in {18,19,20,21,22}{
    \node [scale=0.8, draw, circle, fill=white] at (\x ,1){  } ;
    }
     \foreach \x in {17,18,19,20,21,22}{
    \node [scale=0.8, draw, circle, fill=white] at (\x ,0){  } ;
    }
\end{scriptsize}
\node[below] at (16,-2){(b) :  $f(v_i)=10$,   $a = 13$, $q = 3$} ;
\end{tikzpicture}
\caption{\label{fig:mapping-g-1}Construction of the mapping $g$ in the proof of Lemma~\ref{lem:broadcast-at-most-qa}.}
\end{center}
\end{figure}

We now prove that in all cases, $g$ is an independent broadcast on $\overrightarrow{C}(qa; 1, a)$. For that, we first prove the following claim.

\begin{claim}\label{cl:distance-qa}
For every vertex $v_j$ whose $g$-value is set to~$d$ in Item 1, or set to~$q-1$  in Item 2, we have
$d(v_i,v_j)\le f(v_i)- g(v_j)$.
\end{claim}

\begin{proof}

From {Proposition}~\ref{prop:Distance1}, it is evident that for each vertex whose $g$-value is set to $d$ in Item 1, we have $$d(v_i,v_{j}) = q \le f(v_i) -g(v_{j}).$$
Now, in Item 2 and still using {Proposition}~\ref{prop:Distance1}, we can easily see that the vertex farthest from $v_i$ among the vertices whose $g$-value might be set to $1$ is the vertex $v_j$ with $j= i + q + (p-2)q $, which gives
$$d(v_i,v_{j}) = q +(p-2)q < pq - (q-1)= f(v_i) - g(v_{j}).$$

This concludes the proof of the claim.
\end{proof}
With this assertion, we can conclude that in each item and for every $g$-broadcast vertex $v_j$ modified by the function $g$, the vertices $g$-dominated by this vertex are already $f$-dominated by the vertex $v_i$. Thus, to prove that $g$ is an independent broadcast on $\overrightarrow{C}(qa;1,a)$, it suffices to show that for every two $g$-broadcast vertices $v_i$ and $v_j$ in some item, we have $d(v_{i},v_{j}) >  g(v_{i})$. Furthermore, applying again the distance formula from {Proposition}~\ref{prop:Distance1}, we can easily see that in each item $d(v_i, v_j) \geq q-1 \geq g(v_i)$ for every two $g$-broadcast vertices $v_i$ and $v_j$.

To conclude the proof, we need to prove that we have $\sigma(g) \geq \sigma(f)$. Indeed, in Item~1, the number of vertices set to $q-1$  and to $d$ are respectively $n_1 = 1$ and $n_2 = d+1$, which gives
$$n_1 + n_2(d+1) = q-1 + d(d+1)  \geq f(v_i).$$

In Item~2, the number of vertices set to~$q-1$, are
$n_1 =  1 + q(p-1) $,
which gives
 $$ n_1(q-1) = q(p-1)(q-1) + q-1.$$
Since $ q\geq 3$ and $p\geq 2$, we have $q(p-1)(q-1) + q-1 \geq 2q(p-1) + q - 1 \geq qp + q-1 $  and thus
 $$ n_1(q-1) \geq qp + q-1 = f(v_i).$$

Therefore, we have established $\sigma(g) \geq \sigma(f)$, as required. This concludes the proof.
\end{proof}

Consider now the case of circulant graphs $\overrightarrow{C}(qa+r;1,a)$ with $r=a-1$.
\begin{lemma}\label{lem:broadcast-at-most-q(a-1)-1}
If $n$, $q$ and $a$  are three integers such that $ n= qa +a- 1 $, then $\overrightarrow{C}(n;1,a)$ admits an $(a-1)$-bounded $\beta_b$-broadcast if $a - 1 \le q $, and a $q$-bounded $\beta_b$-broadcast if $a - 1  > q$.
\end{lemma}
\begin{proof}
If $q = 2$, the result follows directly from Theorem~\ref{th:C(3a-1;1,a)}, which states that $\overrightarrow{C}(n;1,a)$ admits a $2$-bounded $\beta_b$-broadcast.
Let $q\geq 3$. We have to prove that, for any independent  broadcast $f$ on $\overrightarrow{G} = \overrightarrow{C}(n;1,a)$, it is always possible to construct another broadcast $g$ on $\overrightarrow{G}$ that satisfies three conditions: Firstly, for every vertex $v \in V_{g}^+$, $g(v)\leq \min \{a-1,q\}$. Secondly, $ g $ is an independent broadcast on $\overrightarrow{G}$, and thirdly,  $ \sigma(g) \geq \sigma(f)$.

Let $f$ be any independent broadcast on $\overrightarrow{G}$. We consider two cases, according to the maximum value between $a$ and $q-1$, to define the mapping $g$ from $V(\overrightarrow{G})$ to $\{0,1 \ldots , \min \{a-1, q\}\}$.

\begin{enumerate}

\item $a - 1 \le q$. In this case, we have $n=q(a+1) - 1 $, and thus the result directly follows from Lemma~\ref{lem:broadcast-at-most-a-1-a<q+r} (The same construction for the function $g$ can be applied here as well).
\item $q < a-1$.
\begin{enumerate}

\item  If $v_i$ is an $f$-broadcast vertex such that $ q+1 \le f(v_{i}) \le 2q-1 $, then we let
$$g(v_j) =
\left\{
   \begin{array}{ll}
   q-1 & \text{if } j=i, \\ [1ex]
  d & \text{if $j= i + (q-k)a + k$ and $ 0 \le k \le d   $, } \\ [1ex]

   \end{array}
\right.
$$
where $d = f(v_i)-q$.

\item  If $v_i$ is an $f$-broadcast vertex such that $pq\le f(v_{i}) \le (p+1)q-1 $, $p \geq 2$,  then we let
\begin{center}
$g(v_j) = q-1 $ if  $j=i$ or ($d  \le j \le d + (p-2)q$ and $(j - d) \equiv 0 \mod (q)$),
\end{center}

where $d = i+(q-k)a+k $ and $ 0 \le k \le q$.

\item For every other vertex $v_i$, we let $g(v_i)=f(v_i)$.
\end{enumerate}
\end{enumerate}

To prove that $g$ is an independent broadcast on $\overrightarrow{C}(qa+a-1; 1, a)$, we proceed in a similar manner as the proof of {Lemma}~\ref{lem:broadcast-at-most-qa}.
\end{proof}

\begin{lemma}\label{lem:vertex l-l-l }
Let $n$, $a$, and $\ell$ be three integers such that $4 \le a$ and $2 \le \ell \le a-1$.
If $\overrightarrow{C}(n;1,a)$ admits an $\ell$-bounded  independent broadcast  $f$, and if there exists some vertex $v_i$ such that $f(v_i) = \ell - 1 $, and $ \ell \le d(v_{i+a-1}, v_i) $, then  the following two properties are satisfied :

\begin{enumerate}
\item If $\ell = a-1$ and  $ 0 \le  f(v_{i + a - 1}) < a-2$,  then there exists an independent broadcast $g$ on $\overrightarrow{C}(n;1,a)$ such that
$$g(v_i) = g(v_{i+a-1}) =  a-2.$$
\item If $\ell < a-1$ and $ 1 \le  f(v_{i + a - 1}) < \ell-1$, then there exists an independent broadcast $g$ on $\overrightarrow{C}(n;1,a)$ such that
 $$g(v_i) = g(v_{i+a-1}) =  \ell - 1.$$
\end{enumerate}
\end{lemma}

\begin{proof}
To prove this lemma, all we need to show is that for any independent broadcast $f$ on $\overrightarrow{C}(n;1,a)$, if there exists a vertex $v_i \in V_f^+$ such that $f(v_i) = \ell-1$ and $\ell \le d(v_{i+a-1}, v_i) $, then there exists an independent broadcast $g$ on $\overrightarrow{C}(n;1,a)$ such that $g(v_i) = g(v_{i+a-1}) = \ell - 1$ and  $\sigma(g)\ge\sigma(f)$. Let $f$ be an $\ell$-bounded independent broadcast on $\overrightarrow{C}(n;1,a)$, where $f(v_i) = \ell - 1$, and  $ \ell \le d(v_{i+a-1}, v_i)$. We define the mapping $g$  from $V(\overrightarrow{C}(n;1,a))$ to  $ \{0,1, \ldots , \ell \}$ as follows.
\begin{enumerate}
\item  $\ell = a-1$.\\
Let $ f(v_{i+a-1})= p $ and $f(v_{i+(p+2)a-1})= q$, such that $ 0 \le p < a-2$ and $ 0 \le q \le a-2$. We consider two subcases, according to  the value of $q$.

\begin{enumerate}
\item If $0 \le q \le a - 3 - p $, then we let
$$g(v_j) =
\left\{
   \begin{array}{ll}
   a-2 & \text{if } j=i, \\ [1ex]
   q + p + 1& \text{if $j=i+a-1$, } \\ [1ex]
    0& \text{if $j=i+(p+2)a-1$. }

   \end{array}
\right.
$$

Note that the vertices $g$-dominated by the vertex $v_{i+a-1}$ are already $f$-dominated by  $v_{i}$, $v_{i+a-1}$ and $v_{i+(p+2)a-1}$, which implies that function $g$ is an independent broadcast on $\overrightarrow{C}(n;1,a)$ of cost $\sigma(g) > \sigma(f)$.

\item If $a-2-p \le q \le a-2 $, then we let
$$g(v_j) =
\left\{
   \begin{array}{ll}
   a-2 & \text{if  $j=i$ or $j=i+a-1$,} \\ [1ex]
   0 & \text{if $j=i+(p+2)a-1$,} \\ [1ex]
    q - (a-2 - p)& \text{if $j=i+(p+2)a + (a-3-p)$.}

   \end{array}
\right.
$$

In this case, the vertices $g$-dominated by the vertex $v_{i+a-1}$ are already $f$-dominated by $v_{i}$, $v_{i+a-1}$, and $v_{i+(p+2)a-1}$. Furthermore, the vertex $v_{i+(p+2)a + (a-3-p)}$ is only $f$-dominated by $v_{i+(p+2)a-1}$, since the vertex $v_{i+(p+1)a + (a-3-p)}$ is $f$-dominated by $v_{i}$. Additionally, since
$$ d(v_{i+(p+2)a-1},v_{i+(p+2)a + (a-3-p)})=  f(v_{i+(p+2)a-1}) - g(i+(p+2)a + (a-3-p)), $$
 then the vertices $g$-dominated by vertex $v_{i+(p+2)a + (a-3-p)}$ are already $f$-dominated by $v_{i+(p+2)a-1}$. This implies that $g$ is an independent broadcast on $\overrightarrow{C}(n;1,a)$, of cost $\sigma(g) = \sigma(f).$

\end{enumerate}

\item  $\ell < a-1$.\\
Let $ f(v_{i+a-1})= p $ and $f(v_{i+(p+2)a-1})= q$, such that $ 0 \le p < \ell -1$ and $ 0 \le q \le \ell-1$. Now, we consider three subcases, depending on the value of $q$

\begin{enumerate}
\item If $0 < p < \ell -1  $ and $0 \le q \le \ell - p - 2  $,  then we let
$$g(v_j) =
\left\{
   \begin{array}{ll}
   \ell -1 & \text{if } j=i, \\ [1ex]
   q + p + 1 & \text{if $j=i+a-1$, } \\ [1ex]
    0& \text{if $j=i+(p+2)a-1$. }

   \end{array}
\right.
$$

Note that the vertices $g$-dominated by the vertex $v_{i+a-1}$ are already $f$-dominated by  $v_{i}$, $v_{i+a-1}$ and $v_{i+(p+2)a-1}$, which implies that $g$ is an independent broadcast on $\overrightarrow{C}(n;1,a)$, of cost $\sigma(g) > \sigma(f)$.

\item If $ 0 < p < \ell -1  $ and $ \ell - p + 1  \le q \le \ell $, then we let
$$g(v_j) =
\left\{
   \begin{array}{ll}
   \ell - 1 & \text{if  $j=i$  or  $j=i+a-1$, } \\ [1ex]
   0 & \text{if $j=i+(p+2)a-1$, } \\ [1ex]
    q - (\ell - 1 - p)& \text{if $j=i+(p+2)a + (a-3-p)$. }
   \end{array}
\right.
$$

The vertices $g$-dominated by the vertex $v_{i+a-1}$ are already $f$-dominated by $v_{i}$, $v_{i+a-1}$, and $v_{i+(p+2)a-1}$. Since the vertex $v_{i+(p+1)a + (a-3-p)}$ is $f$-dominated by $v_{i}$, the vertex $v_{i+(p+2)a + (a-3-p)}$ is only $f$-dominated by $v_{i+(p+2)a-1}$.
Additionally, since
$$ d(v_{i+(p+2)a-1},v_{i+(p+2)a + (a-3-p)})=  f(v_{i+(p+2)a-1}) - g(i+(p+2)a + (a-3-p)),$$
then the vertices $g$-dominated by vertex $v_{i+(p+2)a + (a-3-p)}$ are already $f$-dominated by $v_{i+(p+2)a-1}$, this implies that $g$ is an independent broadcast on $\overrightarrow{C}(n;1,a)$, of cost $\sigma(g) = \sigma(f).$

\item If $ p = 0   $ and $ q \le \ell - 1$, then we distinguish two sub-cases depending on whether the vertex $v_{i+a-1}$ is $f$-dominated or not

\begin{enumerate}
\item If $v_{i+a-1}$ is not $f$-dominated by any vertices, then we let
$$g(v_j) =
\left\{
   \begin{array}{ll}
   \ell & \text{if } j=i, \\ [1ex]
   q+1 & \text{if $j=i+a-1$, } \\ [1ex]
   0 & \text{if $j=i+2a-1$. }

   \end{array}
\right.
$$

In this case, the vertices $g$-dominated by the vertex $v_{i+a-1}$ are already $f$-dominated by  $v_{i}$ and $v_{i+2a-1}$, which implies that function $g$ is an independent broadcast on $\overrightarrow{C}(n;1,a)$, of cost $\sigma(g) > \sigma(f)$.

\item If $v_{i+a-1}$ is $f$-dominated by some vertex $v_k$ with $i < k < i+a $, then we let
$$g(v_j) =
\left\{
   \begin{array}{ll}
   \ell & \text{if } j=i, \\ [1ex]
   f(v_k)-1 & \text{if $j=k$, } \\ [1ex]
   q+1 & \text{if $j=i+a-1$, } \\ [1ex]
   0 & \text{if $j=i+2a-1$. }

   \end{array}
\right.
$$

We can first easily show also that the vertices $g$-dominated by vertex $v_{i+a-1}$ are already $f$-dominated by  $v_{i}$ and $v_{i+2a-1}$. Furthermore, since the vertex $v_{i+a-1}$ is $f$-dominated only by $v_{k}$, we have $d(v_{k},v_{i+a-1})=f(v_{k})$, otherwise, $v_{k}$ would $f$-dominate $v_{i+2a-1}$,
contradicting the fact that $f$ is an independent broadcast. Therefore, the vertex $v_{i+a-1}$ is not $g$-dominated by $v_{k}$. Consequently, we can conclude that $g$ is an independent broadcast on $\overrightarrow{C}(n;1,a)$, of cost $\sigma(g) = \sigma(f).$ Since $g(v_i) = \ell$, $g(v_{i+a-1}) > 0 $, and according to both  Items $(a)$ and $(b)$,  there exists an independent broadcast $h$ on $\overrightarrow{C}(n;1,a)$ such that $h(v_i) = h(v_{i+a-1}) = \ell-1$ and  $\sigma(h)\ge\sigma(f)$.
\end{enumerate}
\end{enumerate}
\end{enumerate}
\end{proof}

\subsection{Bounds and exact values for $\beta_b (\overrightarrow{C}(n;1,a))$}\label{sec:bounds}

In this section, we will present general upper and lower bounds on the cost of independent broadcasts on oriented circulant graph $\overrightarrow{G}=\overrightarrow{C}(n;1,a)$. We will explore these bounds in two distinct contexts: one where the graph $\overrightarrow{G}$ admits an $(a-1)$-bounded independent broadcast, and another where it admits an $q$-bounded independent broadcast.

\subsubsection{Bounds for $(a-1)$-bounded independent broadcast on $\overrightarrow{C}(n;1,a)$}\label{sec:bounds a-1 }

We first introduce some notation and a useful lemma in the case where $ \overrightarrow{G}$ admits an $(a-1)$-bounded independent broadcast. Let $f$ be an $(a-1)$-bounded independent broadcast on $\overrightarrow{G}$. We set
$$V_f^{1}=\{v_i\in V_f\ |\ f(v_i) = a-1 \} \mbox{ and } V_f^{2}=\{v_i\in V_f\ |\ f(v_i) \leq a-2\}.$$
In particular, if $f$ is $(a-1)$-bounded, we have $ V_f^+ = V_f^{1} \cup V_f^{2}$.\\ For each vertex \(v_i\in V_f^1\) and \(v_j\in V_f^2\), let
\begin{center}
$ A_f^i = \{v_{i+k},\ 0\le k \le 2a-2\} $ and $ B_f^j = \{v_{j+k},\ 0\le k \le f(v_j)\}$.
\end{center}

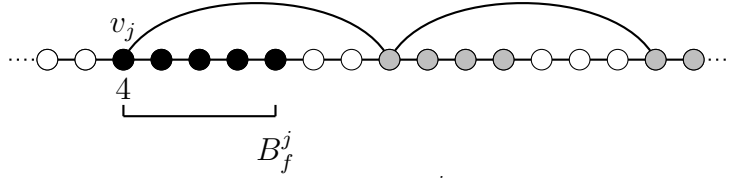
\begin{figure}[h]
\begin{center}
 \begin{tikzpicture}[scale=0.5]

    \POINTILLE{0}{0}{1}{0}   \LIGNE{1}{0}{24}{0}   \POINTILLE{24}{0}{25}{0}

   \foreach \k in {2,9,16} \draw[thick] (\k,0) .. controls (\k+1,2) and (\k+6,2) .. (\k+7,0);
   \foreach \k in {1,2,...,24} \SOM{\k}{0}{}{};

   \bSOM{2}{0}{$v_i$}{6}

   \SOM{1}{0}{}{}
   \bSOM{3}{0}{}{}
   \bSOM{4}{0}{}{}
   \bSOM{5}{0}{}{}
   \bSOM{6}{0}{}{}
   \bSOM{7}{0}{}{}
   \bSOM{8}{0}{}{}
   \bSOM{9}{0}{}{}
   \bSOM{10}{0}{}{}
   \bSOM{11}{0}{}{}
   \bSOM{12}{0}{}{}
   \bSOM{13}{0}{}{}
   \bSOM{14}{0}{}{}
    \gSOM{16}{0}{}{}
   \gSOM{17}{0}{}{}
   \gSOM{18}{0}{}{}
   \gSOM{19}{0}{}{}
   \gSOM{20}{0}{}{}
    \gSOM{23}{0}{}{}
     \gSOM{24}{0}{}{}
    \LIGNE{2}{-1.2}{2}{-1.5}
   \LIGNE{2}{-1.5}{14}{-1.5}
   \LIGNE{14}{-1.5}{14}{-1.2}
 \node[below] at (8,-1.6){ $A_f^i$  } ;
   \node[below] at (12,-2.8){The vertices of the set $A_f^i$, $f(v_i)=6 $, $a=7$ } ;
\end{tikzpicture}
\begin{tikzpicture}[scale=0.5]
   \POINTILLE{0}{0}{1}{0}   \LIGNE{1}{0}{18}{0}   \POINTILLE{18}{0}{19}{0}
   \foreach \k in {3,10} \draw[thick] (\k,0) .. controls (\k+1,2) and (\k+6,2) .. (\k+7,0);
   \foreach \k in {1,2,...,18} \SOM{\k}{0}{}{};
   \SOM{1}{0}{}{}
   \SOM{2}{0}{}{}
   \bSOM{3}{0}{$v_j$}{4}
   \bSOM{4}{0}{}{}
   \bSOM{5}{0}{}{}
   \bSOM{6}{0}{}{}
   \bSOM{7}{0}{}{}
   \gSOM{10}{0}{}{}
   \gSOM{11}{0}{}{}
   \gSOM{12}{0}{}{}
   \gSOM{13}{0}{}{}
   \gSOM{17}{0}{}{}
   \gSOM{18}{0}{}{}
     \LIGNE{3}{-1.2}{3}{-1.5}
   \LIGNE{3}{-1.5}{7}{-1.5}
   \LIGNE{7}{-1.5}{7}{-1.2}
 \node[below] at (7,-1.6){ $B_f^j$  } ;
    \node[below] at (10,-2.8){The vertices of the set $B_f^j$, $f(v_j)=4$, $a=7$ } ;
\end{tikzpicture}

\caption{\label{fig:AiBj}The sets $A_f^i$ and $B_f^j$.}
\end{center}
\end{figure}

The definition of these two sets is illustrated in Figure~\ref{fig:AiBj}.

These sets have the following properties:
\begin{lemma}\label{lem:AiBj}
For every $(a-1)$-bounded independent broadcast $f$ on $\overrightarrow{C}(n;1,a)$, $4\le a$, the following holds.
\begin{enumerate}
\item For every vertex $v_i\in V_f^1$, $|A_f^i| = 2a-1$,
\item For every vertex $v_j\in V_f^2$, $|B_f^j| = f(v_j)+1$,
\item $\sum_{v_i\in V_f^1}|A_f^i| + \sum_{v_j\in V_f^2}|B_f^j|  \le n$.
\end{enumerate}
\end{lemma}

\begin{proof}
The first two items directly follow from the definition of the sets $A_f^i$ and $B_f^j$. It also follows from the definition that $A_f^i \cap A_f^{i'} = \emptyset$ for every
two distinct vertices $v_i$ and $v_{i'}$ in $V_f^1$, since otherwise we would have $d(v_{j},v_{j'})\leq \max \{f(v_{j}),f(v_{j'})\} = a-1$,
contradicting the fact that $f$ is an independent broadcast.
Similarly, we have  $B_f^j \cap B_f^{j'} = \emptyset$
for every two distinct vertices $v_j$ and $v_{j'}$ in $V_f^2$. The same argument gives $A_f^i \cap B_f^{j} = \emptyset$
for every two vertices $v_i\in V_f^1$ and $v_j\in V_f^{2}$.
With these three intersections empty, Item~3 also holds.
This completes the proof.
\end{proof}

\begin{proposition}\label{prop:bound1-a-1}
If $n, a, q$ and $r$ are four integers such that $n=qa+r$. If $4 \le a \leq q+r+1$ and $0 \le r \le a-2$, or $a \le q + 1 $ and $r=a-1$, then, for every $(a-1)$-bounded independent broadcast $f$ on $\overrightarrow{C}(n;1,a)$, we have
$$\sigma(f)\leq  n - \left|V_f^+ \right| - (a-1)\left|V_f^1 \right|.$$
\end{proposition}

\begin{proof}
Let $f$ be any $(a-1)$-bounded independent broadcast on $\overrightarrow{C}(n;1,a)$, which exists from Lemma~\ref{lem:broadcast-at-most-a-1-a<q+r}.
From Lemma~\ref{lem:AiBj}, we get
$$\sum_{v_i\in V_f^1}|A_f^i|  +  \sum_{v_j\in V_f^{ 2}}|B_f^j|
   =  2f(V_f^1) + \left|V_f^1\right| + f(V_f^{ 2}) + \left|V_f^{ 2}\right| \le n,$$
which gives
$$ f(V_f^+) = f(V_f^1) + f(V_f^{ 2}) \le n - \left|V_f^+ \right| - f(V_f^1).$$
Now, since $f(v_i)= a-1$ for every $v_i\in V_f^{ 1}$,
we have $f(V_f^{ 1}) = (a-1)\left|V_f^{ 1}\right|$, and thus
$$\sigma(f) = f(V_f^+) \leq  n - \left|V_f^+ \right| - (a-1) \left|V_f^1 \right|.$$
This completes the proof.
\end{proof}

\begin{proposition}\label{prop:bound2-a-1}
If $n, a, q$ and $r$ are four integers such that $n=qa+r$. If $4 \le a \leq q+r+1$ and $0 \le r \le a-2$, or $a \le q + 1 $ and $r=a-1$,
then, for every $(a-1)$-bounded independent broadcast $f$ on $\overrightarrow{C}(n;1,a)$, we have
$$\sigma(f) \leq \left\lfloor  \left( \frac{a-2}{a-1} \right) \left(n - \left|V_f^{1}\right| \right) \right\rfloor - (a-3)\left|V_f^{1}\right|.$$
\end{proposition}

\begin{proof}
Let $f$ be any $(a-1)$-bounded independent broadcast on $\overrightarrow{C}(n;1,a)$. For every vertex $v_j\in V_f^2$, we have $f(v_j) \leq a-2$ and $|B_f^j| = f(v_j) + 1$, this implies
$$\frac{|B_f^j|}{f(v_j)} = \frac{f(v_j) + 1}{f(v_j)} = 1 + \frac{1}{f(v_j)}
\ge 1 + \frac{1}{a-2} = \frac{a-1}{a-2},$$
which gives
$$|B_f^j| \ge \left( \frac{a-1}{a-2} \right) f(v_j),$$
and thus
$$ \sum_{v_j\in V_f^{ 2}}|B_f^j| \ge \left( \frac{a-1}{a-2} \right)  f(V_f^{ 2}).$$

From Lemma~\ref{lem:AiBj}, we get
$$ n \geq \sum_{v_i\in V_f^1}|A_f^i|  +  \sum_{v_j\in V_f^{ 2}}|B_f^j|
   \geq  2f(V_f^1) + \left|V_f^1\right| + \left( \frac{a-1}{a-2} \right)  f(V_f^{ 2}),$$
which gives
$$n \ge  \left( \frac{a-1}{a-2} \right)f(V)  + \left( \frac{a-3}{a-2} \right )f(V_f^{1}) + \left|V_f^{ 1}\right|.$$

Finally, since $f(V_f^{ 1}) = (a-1)\left|V_f^{1}\right|$, we get

$$\left( \frac{a-1}{a-2} \right)f(V) \leq n - \left|V_f^{ 1}\right| - \left( \frac{a-3}{a-2} \right )(a-1)\left|V_f^{1}\right|,$$

which gives

$$f(V) \leq \left( \frac{a-2}{a-1} \right) \left( n - \left|V_f^{1}\right| \right) -  (a-3)\left|V_f^{1}\right|,$$

and thus

$$\sigma(f) = f(V) \leq \left\lfloor  \left( \frac{a-2}{a-1} \right) \left(n - \left|V_f^{1}\right| \right) \right\rfloor - (a-3)\left|V_f^{1}\right|.$$

This completes the proof.
\end{proof}

\begin{lemma}\label{lem:bound-Inf-a-1}
If $n$, $q$, $k$ and $s$  are four integers such that $ n= k(a-1)+s $, $a\geq 4$ and $0 \le s \le \min\{a,k\}-2$, then $\overrightarrow{C}(n;1,a)$ admits an $(a-1)$-bounded $\beta_b$-broadcast.
\end{lemma}

\begin{proof}
 Let $n=ka - (k-s)$ with $a\geq 4$ and $0 \le  s \le \min\{a,k\}-2$. We consider two cases, depending on the maximum value between $a$ and $k-s$.
\begin{enumerate}
\item $ k-s  \leq a$. We can write $n$ as  $n= qa + r$, where $q=k-1$ and $r=(a-k+s)$. Since $k\ge 3$ and $k-s \ge 2$, we get $q \geq 2$, $0 \le r \le a-2$ and $a\le q+r+1$. From Lemma~\ref{lem:broadcast-at-most-a-1-a<q+r}, $\overrightarrow{C}(n;1,a)$ admits an $(a-1)$-bounded independent broadcast.
\item $ a < k-s$. Let  $ k - s = pa+t$, with $1 \le p$ and $0 \le t \le a-1$, which gives $$ n = k a - (k-s) = (k-p-1)a + a-t.$$
Now, we consider three cases, depending on the value of $t$
\begin{enumerate}
\item If $t=0$, then  $n=(k-p)a $, with $p \geq 2$ because otherwise $k-s=a$, which is impossible. Since $a\geq 4$ then, $ k \ge pa \ge 4p$, and thus   $ k-p \ge 3p \ge 6$, which allows us to denote $n=qa$ with $q = k-p.$ Moreover, since $2 \le p$, $4\le a$ and $k=pa+s$, then, $a \le p(a-1) \le k-p   $, and then  $a \le q$. From Lemma~\ref{lem:broadcast-at-most-qa}, $\overrightarrow{C}(n;1,a)$ admits an $(a-1)$-bounded independent broadcast.

\item If $t=1$, then  $n=(k-p-1)a + a - 1 $, with $p \geq 1.$ Since $a\geq 4$ then, $ k \ge pa + 1 \ge 4p +1 $, and thus   $ k - p -1 \ge 3p \ge 3$, this also allows us to denote $n=qa+r$ with $q = (k-p-1)$ and $r=a-1$. Moreover, since $1 \le p$, $4\le a$, and $k=pa+s$ then,
$a \le p(a-1)+1 \le k-p $, and thus  $a \le q + 1 $. From Lemma~\ref{lem:broadcast-at-most-q(a-1)-1}, $\overrightarrow{C}(n;1,a)$ admits an $(a-1)$-bounded independent broadcast.

\item If $2 \le t \le a-1$, then we can write $n$ as  $n= qa + r$, where $q=k-p-1$, $r=a-t$ and $0 \le r \le a-2$. Since $p\geq 1$ and $a\geq 4$ then $ k \geq pa + 2 \ge 4p+2 \geq p +5  $, and thus $k-p-1 > 3$. Moreover, we have  $q+r+1= k - p + a - t \ge p(a-1) +a \ge a $, then  $a  \le q+r+1$. From Lemma~\ref{lem:broadcast-at-most-a-1-a<q+r}, $\overrightarrow{C}(n;1,a)$ admits an $(a-1)$-bounded independent broadcast.
\end{enumerate}
\end{enumerate}
\end{proof}

\begin{proposition}\label{prop:bound-Inf-a-1}
Let $n$, $a$, $k$, and $s$ be four integers such that $n=k(a-1) + s$ and $a\geq 4$. If  $0 \le s \le \min\{a,k\}-2$, then
\begin{enumerate}
\item $\beta_b(\overrightarrow{C}(n;1,a)) \geq k(a-2), \mbox { for } s=0 $,
\item $\beta_b(\overrightarrow{C}(n;1,a)) \geq (k-1)(a-2)+1, \mbox { for } s=1 $,
\item $\beta_b(\overrightarrow{C}(n;1,a)) \geq (k-s)(a-2)+(s-1)(s+1), \mbox { for } s \geq 2$.
\end{enumerate}
\end{proposition}
\begin{proof}
From Lemma~\ref{lem:bound-Inf-a-1}, $\overrightarrow{C}(n;1,a)$ admits an $(a-1)$-bounded independent broadcast.
We now construct a mapping $f$ from $V(\overrightarrow{C}(n;1,a))$ to $\{0,1, \ldots ,a-1\}$, based on the result of  Lemma~\ref{lem:vertex l-l-l }, by setting the cost $a-2$ to the vertex $v_0$. We consider three cases, depending on the value of $s$.
\begin{enumerate}
\item $s=0$. We let $f(v_i)=a-2$ if $i \equiv 0 \pmod{a-1}$, and $f(v_i)=0$ otherwise. From Proposition~\ref{prop:Distance1}, we have for every two vertices $v_i,v_j \in V_f^+$, $d(v_i,v_j) = a-1$. Therefore, $f$ is clearly an independent broadcast on $\overrightarrow{C}(n;1,a)$.
\item $s=1$. We let
$$f(v_i) =
\left\{
   \begin{array}{ll}
   a-2 & \text{if $i = p(a-1)$ and $0 \le p \le  k-3 $}, \\ [1ex]
   a-1  & \text{if $i=(k-2)(a-1)$}, \\ [1ex]
   0 & \text{otherwise.}
   \end{array}
\right.
$$
From Proposition~\ref{prop:Distance1},  for every two vertices $v_i,v_j \in V_f^+ $ such that $v_i \neq v_{(k-2)(a-1)}$, $d(v_{i},v_j)  = a - 1 > f(v_i) $ and $d( v_{(k-2)(a-1)},v_i)= a > f(v_{(k-2)(a-1)})$. Therefore, $f$ is an independent broadcast on $\overrightarrow{C}(n;1,a)$.
\item $2\leq s \leq a-1$. We let
$$f(v_i) =
\left\{
   \begin{array}{ll}
   a-2 & \text{if $i = p(a-1)$ and $0 \le p \le  k - s - 1 $}, \\ [1ex]
   s-1  & \text{if $i = p(a-1)$ and $k - s  \le p  \le k $}, \\ [1ex]
   0 & \text{otherwise.}
   \end{array}
\right.
$$	
From Proposition~\ref{prop:Distance1}, for every two vertices $v_i,v_j \in V_f^+ $ such that $ 0 \le i \le (k-s-1)(a-1)$, $d(v_{i},v_j)  = a - 1 > f(v_i) $, and for every two vertices $v_i,v_j \in V_f^+ $ such that $ (k-s)(a-1) \le i \le k(a-1)$,  $d( v_i,v_j)= s-1 > f(v_{i})$. Therefore, $f$ is an independent broadcast on $\overrightarrow{C}(n;1,a)$
\end{enumerate}
\end{proof}

\subsubsection{Bounds for optimal $q$-bounded independent broadcast on $\overrightarrow{C}(n;1,a)$}\label{sec: bounds q}
In this section, we consider oriented circulant graphs $\overrightarrow{G}=\overrightarrow{C}(n;1,a)$ such that $n=qa$ or $n=qa+a-1$,  with $3 \le q+1<a$. From Lemma~\ref{lem:broadcast-at-most-qa} and Lemma~\ref{lem:broadcast-at-most-q(a-1)-1},  $\overrightarrow{G}$ admits an $q$-bounded independent broadcast, say   $f$. We set
 $$V_f^{1}=\{v_i\in V_f\ |\ f(v_i) = q \} \mbox{ and } V_f^{2}=\{v_i\in V_f\ |\ f(v_i) \leq q-1\}.$$ In particular, if $f$ is $q$-bounded, we then have $ V_f^+ = V_f^{1} \cup V_f^{2}$. \\
 We distinguish two cases depending on the value of $n$.
\begin{enumerate}
\item If $n=qa$, then, for each vertex \(v_i\in V_f^1\) and \(v_j\in V_f^2\), we set
\begin{center}
$
A_f^i = \{v_{i+k},$ with $\ 0\le k \le q $ or ($k= i+pa +q-p $ and $ 1 \le p \le q-1) \}$ and
\end{center}
\begin{center}
 $ B_f^j = \{v_{j+k},\ 0\le k \le f(v_j)\}$.
\end{center}
\item If $n=qa+a-1$ then, for each vertex \(v_i\in V_f^1\) and \(v_j\in V_f^2\), we set
\begin{center}
$
A_f^i = \{v_{i+k},$ with $\ 0\le k \le q $ or ($k= i+pa +q-p $ and $ 1 \le p \le q) \}$ and
\end{center}
\begin{center}
 $ B_f^j = \{v_{j+k},\ 0\le k \le f(v_j)\}$.
\end{center}
\end{enumerate}

The definition of these two sets is illustrated in Figure~\ref{fig:AiBjq}.
These sets have the following properties.

\begin{figure}[h]
\begin{center}
\begin{tikzpicture}[scale=0.4]
\draw(0,0)circle(5) ;
\LIGNE{0}{5}{4}{-3}
\LIGNE{4}{-3}{-4.75}{-1.5}
   \POINTILLE{-4.}{-3}{0}{5}
   \POINTILLE{4.75}{-1.5}{-4.}{-3}
   \bSOM{0}{5}{2}{}
   \bSOM{2}{4.5}{}{}
   \bSOM{3.5}{3.5}{}{}
   \SOM{4.5}{2}{}{}  \node[below] at (5.2,2.8) { };
   \SOM{5}{0.3}{}{}
   \SOM{4.75}{-1.5}{}{} \node[below] at (5.5,-0.8) { };
   \gSOM{4}{-3}{}{}
   \bSOM{2.7}{-4.2}{}{}
   \SOM{1}{-4.8}{}{ }
  \SOM{-1}{-4.8}{}{}
  \SOM{-2.7}{-4.2}{}{ }
  \SOM{-4}{-3}{}{}
   \bSOM{-4.75}{-1.5}{}{}
   \SOM{-5}{0.3}{}{} \node[below] at (-5.8,0.8) {1};
    \SOM{-4.5}{2}{}{}
  \SOM{-3.5}{3.5}{}{} \node[below] at (-4.1,4.5) {1};
     \SOM{-2}{4.5}{}{}

  \node[rotate=-58] at (3.74,-2.5) {\textbf{{\large $>$}}};
  \node[rotate=170] at (-4.2,-1.58) {\textbf{{\large $>$}}};
\node[below] at (0,-6){The vertices of the set $A_f^i$   } ;
\node[below] at (0,-7.5){with  $n=qa+q-1$,};
\node[below] at (0,-9){ $f(v_i)=2 $, $q=2$ and $a=6$ } ;
\draw(14,0)circle(5) ;
\LIGNE{14}{5}{18.5}{-2}
\LIGNE{18.5}{-2}{9.5}{-2}
\LIGNE{14}{5}{9.5}{-2}
   \bSOM{14}{5}{3}{}
   \bSOM{16}{4.5}{}{}
   \bSOM{17.5}{3.5}{}{}
   \bSOM{18.5}{2}{}{}  \node[below] at (5.2,2.8) { };
   \SOM{19}{0.3}{}{}
    \gSOM{18.5}{-2}{}{}
     \gSOM{17.2}{-3.8}{}{}
     \SOM{10.8}{-3.8}{}{}
     \bSOM{15}{-4.8}{}{}
     \SOM{12.6}{-4.8}{}{}
   \gSOM{9.5}{-2}{}{}
   \SOM{12}{4.5}{}{}
   \SOM{10.5}{3.5}{}{}
   \SOM{9.5}{2}{}{}  \node[below] at (5.2,2.8) { };
   \bSOM{9}{0.3}{}{}
   \node[below] at (14,-6){The vertices of the set $A_f^i$  } ;
\node[below] at (14,-7.5){with $n=qa$, $f(v_i)=3 $, } ;
 \node[below] at (14,-9){ $q=3$ and $a=5$ } ;
\draw(28,0)circle(5) ;
\LIGNE{28}{5}{32.5}{-2}
\LIGNE{32.5}{-2}{23.5}{-2}
\LIGNE{28}{5}{23.5}{-2}
   \bSOM{28}{5}{2}{}
   \bSOM{30}{4.5}{}{}
   \bSOM{31.5}{3.5}{}{}
   \SOM{32.5}{2}{}{}  \node[below] at (19.2,2.8) { };
   \SOM{33}{0.3}{}{}
    \gSOM{32.5}{-2}{}{}
     \gSOM{31.2}{-3.8}{}{}
     \SOM{24.8}{-3.8}{}{}
     \SOM{29}{-4.8}{}{}
     \SOM{26.6}{-4.8}{}{}
   \gSOM{23.5}{-2}{}{}
   \SOM{26}{4.5}{}{}
   \SOM{24.5}{3.5}{}{}
   \SOM{23.5}{2}{}{}  \node[below] at (19.2,2.8) { };
   \SOM{23}{0.3}{}{}
   \node[below] at (28,-6){The vertices of the set $B_f^j$  } ;
\node[below] at (28,-7.5){with $n=qa$, $f(v_i)=2 $ } ;
 \node[below] at (28,-9){ $q=3$ and $a=5$} ;
\end{tikzpicture}
\caption{\label{fig:AiBjq} The sets $A_f^i$ and $B_f^j$ .}
\end{center}
\end{figure}

\begin{lemma}\label{lem:AiBjq}
For every $q$-bounded independent broadcast $f$ on $\overrightarrow{C}(n;1,a)$, the following holds.
\begin{enumerate}
\item For every vertex $v_i\in V_f^1$ we have

$$ |A_f^i| =
\left\{
   \begin{array}{ll}
   2q & \text{if $n = qa$,  } \\ [1ex]
   2q+1  & \text{if $n=qa+a-1$}.
   \end{array}
\right.
$$	
\item For every vertex $v_j\in V_f^2$, $|B_f^j| = f(v_j)+1$.
\item $\sum_{v_i\in V_f^1}|A_f^i| + \sum_{v_j\in V_f^2}|B_f^j|  \le n$.
\end{enumerate}
\end{lemma}

\begin{proof}
The first two items directly follow from the definition of the sets $A_f^i$ and $B_f^j$.
It also follows from the definition that $A_f^i \cap A_f^{i'} = \emptyset$ for every
two distinct vertices $v_i$ and $v_{i'}$ in $V_f^1$. Similarly,  from the definition we have  $B_f^j \cap B_f^{j'} = \emptyset$
for every two distinct vertices $v_j$ and $v_{j'}$ in $V_f^2$,
The same argument gives $A_f^i \cap B_f^{j} = \emptyset$
for every two vertices $v_i\in V_f^1$ and $v_j\in V_f^{2}$.
All together, these three properties imply that Item~3 also holds.
\end{proof}

\begin{proposition}\label{prop:bound1-q}
If $n$, $a$, $q$ and $r$ are four integers such that $n=qa+r$, $r \in \{0,a-1\} $ and $q < a-1 $ then, for every independent broadcast $f$ on $\overrightarrow{C}(n;1,a)$, we have

$$ \sigma(f) = f(V_f^+) \leq
\left\{
   \begin{array}{ll}
   n - \left|V_f^2 \right| - q \left|V_f^1 \right| & \text{if $r = 0$,  }  \\ [1ex]
   n - \left|V_f^+ \right| - q \left|V_f^1 \right|  & \text{if $r=a-1$}.
   \end{array}
\right.
$$	

$$ .$$
\end{proposition}

\begin{proof}
Let $f$ be any $q$-bounded independent broadcast on $\overrightarrow{C}(n;1,a)$, which exists from Lemma~\ref{lem:broadcast-at-most-q(a-1)-1} and Lemma~\ref{lem:broadcast-at-most-qa}. We consider two cases, depending on the value of $r$.
\begin{enumerate}
\item $r=0$. From Lemma~\ref{lem:AiBjq}, we get
$$\sum_{v_i\in V_f^1}|A_f^i|  +  \sum_{v_j\in V_f^{ 2}}|B_f^j|
   =  2f(V_f^1) + f(V_f^{ 2}) + \left|V_f^{ 2}\right| \le n,$$
which gives
$$ f(V_f^+) = f(V_f^1) + f(V_f^{ 2}) \le n - \left|V_f^2 \right| - f(V_f^1).$$
Now, since $f(v_i)= q$ for every $v_i\in V_f^{ 1}$,
we have $f(V_f^{ 1}) = q\left|V_f^{ 1}\right|$, and thus
$$\sigma(f) = f(V_f^+) \leq  n - \left|V_f^2 \right| - q \left|V_f^1 \right|.$$
\item$r=a-1$. From Lemma~\ref{lem:AiBjq}, we get
$$\sum_{v_i\in V_f^1}|A_f^i|  +  \sum_{v_j\in V_f^{ 2}}|B_f^j|
   =  2f(V_f^1) + \left|V_f^{1}\right| + f(V_f^{ 2}) + \left|V_f^{ 2}\right| \le n,$$
which gives
$$ f(V_f^+) = f(V_f^1) + f(V_f^{ 2}) \le n - \left|V_f \right| - f(V_f^1).$$
Now, since $f(v_i)= q$ for every $v_i\in V_f^{ 1}$,
we have $f(V_f^{ 1}) = q\left|V_f^{ 1}\right|$, and thus
$$\sigma(f) = f(V_f^+) \leq  n - \left|V_f^+ \right| - q \left|V_f^1 \right|.$$
\end{enumerate}

This completes the proof.
\end{proof}

\begin{proposition}\label{prop:bound2-q}
If $n$, $a$, $q$ and $r$ are four integers such that $n=qa+r$, $r \in \{0,a-1\} $ and $q < a-1 $ then, for every independent broadcast $f$ on $\overrightarrow{C}(n;1,a)$, we have
$$ \sigma(f) = f(V_f^+) \leq
\left\{
   \begin{array}{ll}
   \left\lfloor  \left( \dfrac{q-1}{q} \right) n \right\rfloor - (q-2)\left|V_f^{1}\right|, & \text{if $r = 0$,   } \\ [4ex]
   \left\lfloor  \left( \dfrac{q-1}{q} \right) \left(n - \left|V_f^{1}\right| \right) \right\rfloor - (q-2)\left|V_f^{1}\right|,  & \text{if $r=a-1$}.
   \end{array}
\right.
$$	
\end{proposition}

\begin{proof}
Let $f$ be any $q$-bounded independent broadcast on $\overrightarrow{C}(n;1,a)$, which exists from Lemma~\ref{lem:broadcast-at-most-q(a-1)-1} and Lemma~\ref{lem:broadcast-at-most-qa}. For every vertex $v_j\in V_f^2$, we have $f(v_j) \leq q-1$ and $|B_f^j| = f(v_j) + 1$. Then,

$$\frac{|B_f^j|}{f(v_j)} = \dfrac{f(v_j) + 1}{f(v_j)} = 1 + \frac{1}{f(v_j)}
\ge 1 + \frac{1}{q-1} = \dfrac{q}{q-1},$$
which gives
$$|B_f^j| \ge \left( \dfrac{q}{q-1} \right) f(v_j),$$
and thus
$$ \sum_{v_j\in V_f^{ 2}}|B_f^j| \ge \left( \dfrac{q}{q-1} \right)  f(V_f^{ 2}).$$

We consider two cases, depending on the value of $r$.
\begin{enumerate}
\item $r=0$. From Lemma~\ref{lem:AiBjq}, we get
$$ n \geq \sum_{v_i\in V_f^1}|A_f^i|  +  \sum_{v_j\in V_f^{ 2}}|B_f^j|
   \geq  2f(V_f^1) + \left(\dfrac{q}{q-1} \right)  f(V_f^{ 2}),$$
which gives

$$n \ge  \left( \dfrac{q}{q-1} \right)f(V)  + \left( \dfrac{q-2}{q-1} \right )f(V_f^{1}).$$

Since $f(V_f^{ 1}) = q\left|V_f^{1}\right|$, we get

$$\left( \dfrac{q}{q-1} \right)f(V) \leq n - \left( \dfrac{q-2}{q-1} \right )q\left|V_f^{1}\right|,  $$

which gives

$$f(V) \leq \left( \dfrac{q-1}{q} \right)  n  -  (q-2)\left|V_f^{1}\right|, $$

and thus

$$\sigma(f) = f(V) \leq \left\lfloor  \left( \dfrac{q-1}{q} \right) n \right\rfloor - (q-2)\left|V_f^{1}\right|.$$

\item $r=a-1$. From Lemma~\ref{lem:AiBjq}, we get
$$ n \geq \sum_{v_i\in V_f^1}|A_f^i|  +  \sum_{v_j\in V_f^{ 2}}|B_f^j|
   \geq  2f(V_f^1) + \left|V_f^1\right| + \left(\frac{q}{q-1} \right)  f(V_f^{ 2}),$$
which gives

$$n \ge  \left( \frac{q}{q-1} \right)f(V)  + \left( \frac{q-2}{q-1} \right )f(V_f^{1}) + \left|V_f^{ 1}\right|.$$

Finally, since $f(V_f^{ 1}) = q\left|V_f^{1}\right|$, we get

$$\left( \frac{q}{q-1} \right)f(V) \leq n - \left|V_f^{ 1}\right| - \left( \frac{q-2}{q-1} \right )q\left|V_f^{1}\right|,$$

which gives

$$f(V) \leq \left( \frac{q-1}{q} \right) \left( n - \left|V_f^{1}\right| \right) -  (q-2)\left|V_f^{1}\right|,$$

and thus

$$\sigma(f) = f(V) \leq \left\lfloor  \left( \frac{q-1}{q} \right) \left(n - \left|V_f^{1}\right| \right) \right\rfloor - (q-2)\left|V_f^{1}\right|.$$
\end{enumerate}

This completes the proof.
\end{proof}

\begin{proposition}\label{prop:bound-Inf-qa}

Let $n$, $a$, $q$ be three integers with $n=qa$ and $3 \le q < a-1$. If  $a=kq+1+s$ with $k \geq 1 $ and $0 \le s \le q-1$, then we have

\begin{enumerate}
\item $\beta_b(\overrightarrow{C}(n;1,a)) \geq a(q-1), \mbox { for } s=0 $,
\item $\beta_b(\overrightarrow{C}(n;1,a)) \geq (a-1)(q-1)+1, \mbox { for } s=1 $,
\item $\beta_b(\overrightarrow{C}(n;1,a)) \geq (a-s)(q- 1)+s(s-1), \mbox { for } s \geq 2$.
\end{enumerate}
\end{proposition}
\begin{proof}
From Lemma~\ref{lem:broadcast-at-most-qa} and since $q+1 < a$, $\overrightarrow{C}(qa;1,a)$ admit a $q$-bounded independent broadcast.
Thanks Lemma~\ref{lem:vertex l-l-l }, we construct a mapping $f$ from $V(\overrightarrow{C}(qa;1,a))$ to $\{0,1, \ldots,q\}$. We consider three cases, depending on the value of $s$.

\begin{enumerate}
\item $s=0$. We let $f(v_i)=q-1$ if and only if $i \equiv 0 \pmod{q}$ (see Figure~\ref{fig:q cas 1} (a) for the case $a = 7$, $q = 3$ and $k = 2$).

\begin{figure}
\begin{center}
\begin{tikzpicture}[scale=0.6]
\begin{scriptsize}
 \foreach \i in {0,1,...,6} {
        \draw  (\i,-1.5) -- (\i,3);
 }
   \foreach \i in {0,1,2} {
        \draw  (-1,\i) -- (7.5,\i);
 }
 \foreach \x in {-1,...,7}{
 \node [scale=0.8, draw, circle, fill=lightgray] at (\x , 0){  } ;
 \node [scale=0.8, draw, circle, fill=lightgray] at (\x  , 1){  } ;
 \node [scale=0.8, draw, circle, fill=lightgray] at (\x  , 2){  } ;
 }
 \foreach \x in {0,...,6}{
 \node [scale=0.8, draw, circle, fill=lightgray] at (\x , 3){  } ;
 \node [scale=0.8, draw, circle, fill=lightgray]  at (\x , -1){  } ;
 }
\foreach \x in {0,...,6}{
 \node [scale=0.8, draw, circle, fill=white]  at (\x , 3){  } ;
\node [scale=0.8, draw, circle, fill=white]  at (\x , -1){  } ;
}
\foreach \x in {0,1,2}{
 \node [scale=0.8, draw, circle, fill=white]  at (-1,\x ){  } ;
\node [scale=0.8, draw, circle, fill=white]  at (7,\x){  } ;
}
  \foreach \x in {1,4}{
 \node [scale=0.8, draw, circle, fill=black] at (\x , 0){  } ;
 \node at (\x -0.5,0.5){ $2$ } ;
 }
   \foreach \x in {2,5}{
 \node [scale=0.8, draw, circle, fill=black] at (\x , 1){  } ;
 \node at (\x -0.5,1.5){ $2$ } ;
 }
 \foreach \x in {0,3,6}{
 \node [scale=0.8, draw, circle, fill=black] at (\x , 2){  } ;
 \node at (\x -0.5,2.5){ $2$ } ;
 }
   \foreach \x in {0,1,...,6} {
   \node [rotate=-90] at (\x ,-1.5){ $>$ } ;
}
   \foreach \x in {0,1,2} {
   \node  at ( 7.5,\x){ $>$ } ;
}
\end{scriptsize}

\node[below] at (3,-2){(a) : $a = 7$, $q = 3$, $k=2$, $s=0$} ;
\begin{scriptsize}
\foreach \i in {10,...,20} {
        \draw  (\i,-1.5) -- (\i,3);
 }
   \foreach \i in {0,1,2} {
        \draw  (9,\i) -- (21.5,\i);
 }
 \foreach \x in {9,...,21}{
 \node [scale=0.8, draw, circle, fill=lightgray] at (\x , 0){  } ;
 \node [scale=0.8, draw, circle, fill=lightgray] at (\x  , 1){  } ;
 \node [scale=0.8, draw, circle, fill=lightgray] at (\x  , 2){  } ;
 }
 \foreach \x in {10,...,20}{
 \node [scale=0.8, draw, circle, fill=lightgray] at (\x , 3){  } ;
 \node [scale=0.8, draw, circle, fill=lightgray]  at (\x , -1){  } ;
 }
\foreach \x in {10,...,20}{
 \node [scale=0.8, draw, circle, fill=white]  at (\x , 3){  } ;
\node [scale=0.8, draw, circle, fill=white]  at (\x , -1){  } ;
}
\foreach \x in {0,1,2}{
 \node [scale=0.8, draw, circle, fill=white]  at (9,\x ){  } ;
\node [scale=0.8, draw, circle, fill=white]  at (21,\x){  } ;
}
  \foreach \x in {12,15,18}{
 \node [scale=0.8, draw, circle, fill=black] at (\x , 0){  } ;
 \node at (\x -0.5,0.5){ $2$ } ;
 }
   \foreach \x in {13,16,19}{
 \node [scale=0.8, draw, circle, fill=black] at (\x , 1){  } ;
 \node at (\x -0.5,1.5){ $2$ } ;
 }
 \foreach \x in {14,17,20}{
 \node [scale=0.8, draw, circle, fill=black] at (\x , 2){  } ;
 \node at (\x -0.5,2.5){ $2$ } ;
 }
 \node [scale=0.8, draw, circle, fill=black] at (10 , 2){  } ;
 \node at (9.5,2.5){ $3$ } ;
   \foreach \x in {10,11,...,20} {
   \node [rotate=-90] at (\x ,-1.5){ $>$ } ;
}
   \foreach \x in {0,1,2} {
   \node  at ( 21.5,\x){ $>$ } ;
}
\end{scriptsize}
\node[below] at (15,-2){(b) : $a = 11$, $q = 3$, $k=3$, $s=1$} ;
\node[below] at (10,-3.5){ } ;
\end{tikzpicture}

\begin{tikzpicture}[scale=0.6]
\begin{scriptsize}
 \foreach \i in {0,1,...,15} {
        \draw  (\i,-1.5) -- (\i,4);
 }
   \foreach \i in {0,1,2,3} {
        \draw  (-1,\i) -- (16.5,\i);
 }
 \foreach \x in {-1,...,16}{
 \node [scale=0.8, draw, circle, fill=lightgray] at (\x , 0){  } ;
 \node [scale=0.8, draw, circle, fill=lightgray] at (\x  , 1){  } ;
 \node [scale=0.8, draw, circle, fill=lightgray] at (\x  , 2){  } ;
  \node [scale=0.8, draw, circle, fill=lightgray] at (\x  , 3){  } ;
 }
 \foreach \x in {0,...,15}{
 \node [scale=0.8, draw, circle, fill=lightgray] at (\x , 4){  } ;
 \node [scale=0.8, draw, circle, fill=lightgray]  at (\x , -1){  } ;
 }
\foreach \x in {0,...,15}{
 \node [scale=0.8, draw, circle, fill=white]  at (\x , 4){  } ;
\node [scale=0.8, draw, circle, fill=white]  at (\x , -1){  } ;
}
\foreach \x in {0,1,2,3}{
 \node [scale=0.8, draw, circle, fill=white]  at (-1,\x ){  } ;
\node [scale=0.8, draw, circle, fill=white]  at (16,\x){  } ;
}
  \foreach \x in {4,8,12}{
 \node [scale=0.8, draw, circle, fill=black] at (\x , 0){  } ;
 \node at (\x -0.5,0.5){ $3$ } ;
 }
   \foreach \x in {5,9,13}{
 \node [scale=0.8, draw, circle, fill=black] at (\x , 1){  } ;
 \node at (\x -0.5,1.5){ $3$ } ;
 }
 \foreach \x in {6,10,14}{
 \node [scale=0.8, draw, circle, fill=black] at (\x , 2){  } ;
 \node at (\x -0.5,2.5){ $3$ } ;
 }
 \foreach \x in {0,7,11,15}{
 \node [scale=0.8, draw, circle, fill=black] at (\x , 3){  } ;
 \node at (\x -0.5,3.5){ $3$ } ;
 }
   \node [scale=0.8, draw, circle, fill=black] at (1,0){  } ;
    \node [scale=0.8, draw, circle, fill=black] at (2,1){  } ;
     \node [scale=0.8, draw, circle, fill=black] at (3,2){  } ;
    \node at (0.5,0.5){ $2$ }  ;
  \node at (1.5,1.5){ $2$ }  ;
  \node at (2.5,2.5){ $2$ }  ;
   \foreach \x in {0,1,...,15} {
   \node [rotate=-90] at (\x ,-1.5){ $>$ } ;
}
   \foreach \x in {0,1,2,3} {
   \node  at ( 16.5,\x){ $>$ } ;
}
\end{scriptsize}
\node[below] at (7.5,-2){(c) : $a = 16$, $q = 4$, $k=3$, $s=3$} ;
\end{tikzpicture}
\caption{\label{fig:q cas 1}   Construction of $f$
in the proof of Proposition~\ref{prop:bound-Inf-qa}}
\end{center}
\end{figure}
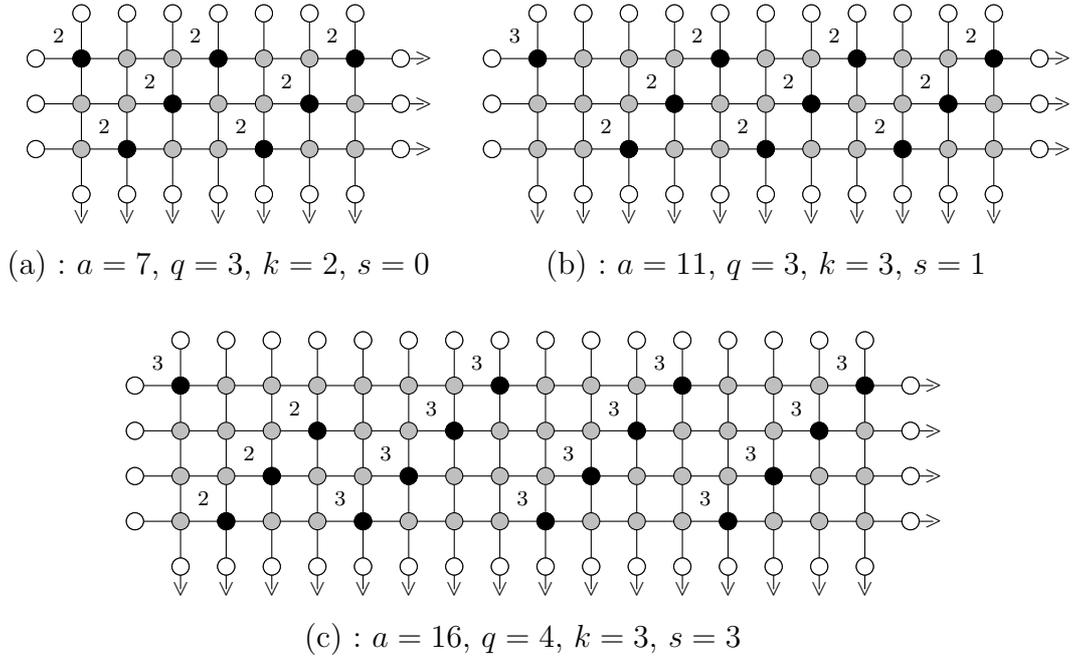

From Proposition~\ref{prop:Distance1},   for every two $g$-broadcast vertices $v_i  \neq v_j$, we have $d(v_i,v_j)=q$. Then, $f$ is an independent broadcast on $\overrightarrow{C}(qa;1,a)$, with cost
$$\sigma(f) =\sum_{v\in V_f^+} f(v) = a(q-1) \le \beta_b(\overrightarrow{C}(n;1,a)).$$

\item$s=1$. We let
$$f(v_i) =
\left\{
   \begin{array}{ll}
    q & \text{if $i=0$, }  \\ [1ex]
      q-1 & \text{if  $i= \alpha(a-1) - \beta q $ with $ 1 \le \alpha \le q$ and $0 \le \beta \le k-1$},\\ [1ex]
     0 & \text{otherwise}.
   \end{array}
\right.
$$
(See Figure~\ref{fig:q cas 1} (b) for the case $a = 11$, $q = 3$ and $k = 3$)

From  Proposition~\ref{prop:Distance1},  for every two $g$-broadcast vertices $v_i \neq v_j$ we have $d(v_i,v_j) \geq f(v_i)+1$. Then $f$ is an independent broadcast on $\overrightarrow{C}(qa;1,a)$, with cost
$$\sigma(f) =\sum_{v\in V_f^+} f(v) = q + kq(q-1) = (a-1)(q-1) + 1 \le \beta_b(\overrightarrow{C}(n;1,a)).$$

\item  $2 \le s \le q-1$. We let
$$f(v_i) =
\left\{
   \begin{array}{ll}
      q-1 & \text{if $i=0$ or  $i= \alpha(a-1) - \beta q $ with $ 1 \le \alpha \le q$ and $0 \le \beta \le k-1$,} \\ [1ex]
  s-1 & \text{if  $i= (q-p)a+p $ with $ 1 \le p \le s $,} \\ [1ex]
   0 & \text{otherwise}.
   \end{array}
\right.
$$
(See Figure~\ref{fig:q cas 1} (c) for the case $a = 16$, $q = 4$ $k=3$ and $s = 3$)

From  Proposition~\ref{prop:Distance1}, we have for every two $g$-broadcast vertices $v_i \neq v_j$ we have $d(v_i,v_j) \geq f(v_i)+1$, this implies that $f$ is an independent broadcast on $\overrightarrow{C}(qa;1,a)$, with cost
$$\sigma(f) =\sum_{v\in V_f^+} f(v) = (kq+1)(q-1) +s(s-1) = (a-s)(q-1) + s(s-1) \le \beta_b(G).$$
\end{enumerate}
This completes the proof.
\end{proof}

\begin{proposition}\label{prop:bound-Inf-qa-1}
Let $n$, $a$, $q$ are three integers with $n=qa+a-1$, $3 \le q < a-1$. If $a=kq+1+s$ with $1 \le k$ and $0 \le s \le q-1$, then
\begin{enumerate}
\item $\beta_b(\overrightarrow{C}(n;1,a)) \geq (a+k)(q-1), \mbox { for } s=0 $,
\item $\beta_b(\overrightarrow{C}(n;1,a)) \geq (a+k-1)(q-1)+1, \mbox { for } s=1 $,
\item $\beta_b(\overrightarrow{C}(n;1,a)) \geq (a+k-s)(q-1)+(s-1)(s+1), \mbox { for } s \geq 2$.
\end{enumerate}
\end{proposition}
\begin{proof} From Lemma~\ref{lem:broadcast-at-most-q(a-1)-1} and since $q+1 < a$, $\overrightarrow{C}(qa +a - 1;1,a)$ admits a $q$-bounded independent broadcast. Thanks to Lemma~\ref{lem:vertex l-l-l }, we construct a mapping $f$ from $V(\overrightarrow{C}(qa+a-1;1,a))$ to $\{0,1, \ldots,q\}$. For this, we consider three cases, depending on the value of $s$.

\begin{enumerate}
\item $s=0$. We let $f(v_i)=q-1$ if and only if $i \equiv 0 \pmod{q}$.
From Proposition~\ref{prop:Distance1}, for every two $g$-broadcast vertices $v_i  \neq v_j$, $d(v_i,v_j)=q$. This implies that $f$ is an independent broadcast on $\overrightarrow{C}(n;1,a)$, with cost
$$\sigma(f) =\sum_{v\in V_f^+} f(v) = (qa+kq)\frac{q-1}{q}=(a+k)(q-1) \le \beta_b(\overrightarrow{C}(n;1,a)).$$

\item  $s=1$. We let
$$f(v_i) =
\left\{
   \begin{array}{ll}
    q & \text{if $i=0$, }  \\ [1ex]
      q-1 & \text{if  $i= \alpha(a-1) - \beta q $ with $ 1 \le \alpha \le q + 1$ and $0 \le \beta \le k-1$}, \\ [1ex]

   0 & \text{otherwise}.

   \end{array}
\right.
$$

From Proposition~\ref{prop:Distance1}, for every two $g$-broadcast vertices $v_i \neq v_j$,   $d(v_i,v_j) \geq f(v_i)+1$. This implies that $f$ is an independent broadcast on $\overrightarrow{C}(qa;1,a)$, with cost
$$\sigma(f) =\sum_{v\in V_f^+} f(v) = q + k(q+1)(q-1) = (a+k-1)(q-1) + 1 \le \beta_b(\overrightarrow{C}(n;1,a)).$$

\item $2 \le s \le q-1$. We let
$$f(v_i) =
\left\{
   \begin{array}{ll}
    q-1 & \text{if $i=0$ or $i= \alpha(a-1) - \beta q $ with $ 1 \le \alpha \le q + 1$ and $0 \le \beta \le k-1$}, \\ [1ex]
  s-1 & \text{if  $i= (q-p)a+ p $ with $0 \le p \le s $}, \\ [1ex]
   0 & \text{otherwise}.

   \end{array}
\right.
$$

From  Proposition~\ref{prop:Distance1},  for every two $g$-broadcast vertices $v_i \neq v_j$, $d(v_i,v_j) \geq f(v_i) +1 $.  This implies that $f$ is an independent broadcast on $\overrightarrow{C}(qa;1,a)$, with cost
\begin{align*}
\sigma(f)
  & =\sum_{v\in V_f^+} f(v) = (k(q+1)+1)(q-1) + (s+1)(s-1)
  \\[1ex]
  & =(a+k-s)(q-1) + (s+1)(s-1) \le \beta_b(\overrightarrow{C}(n;1,a)).
 \end{align*}

\end{enumerate}
This completes the proof.
\end{proof}

\subsubsection{Some exact values}\label{sec:exact}

\begin{theorem}\label{th:ex 1act 1}
If $n$, $a$ and $k$ are three integers such that $n=k(a-1)$, $a\ge 4$ and $k\geq 3$, then
$$\beta_b(\overrightarrow{C}(n;1,a)) = \beta_b(C(n;1,-(n-a)) = (a-2)k. $$

\end{theorem}
\begin{proof}
From Proposition~\ref{prop:bound-Inf-a-1}, $\overrightarrow{C}((a-1)k;1,a)$ admits an $(a-1)$-bounded independent broadcast, which implies
 $$\beta_b(\overrightarrow{C}(n;1,a)) \geq (a-2)k.$$
Moreover, from Proposition~\ref{prop:bound2-a-1}, we get that

$$\sigma(f) \leq \left\lfloor  \left( \frac{a-2}{a-1} \right) \left(n - \left|V_f^{1}\right| \right) \right\rfloor - (a-3)\left|V_f^{1}\right| \leq \left\lfloor  \left( \frac{a-2}{a-1} \right)  n \right\rfloor $$
for every $(a-1)$-bounded independent broadcast $f$ on $\overrightarrow{C}(n;1,a)$,
which gives
$$\beta_b(\overrightarrow{C}(n;1,a)) \le (a-2)k.$$
Finally, by Observation~\ref{isomorphisme}, we then get
$$\beta_b(\overrightarrow{C}(n;1,a)) = \beta_b(\overrightarrow{C}(n;1,-n+a)) =  (a-2)k.$$
This completes the proof.
\end{proof}

\begin{theorem}

If $n$, $a$, $q$ and $k$ are four integers such that $n=qa$, $3 \le q$ and $4 \le a$, then

$$ \beta_b(\overrightarrow{C}(n;1,a))  = \beta_b(\overrightarrow{C}(n;1,-(q-1)a))  =
\left\{
   \begin{array}{ll}
      (q + k)(a-2) & \text{if $q=k(a-1)$, $k\geq 1$, } \\ [1ex]
      a(q-1) & \text{if $a=kq+1$, $k\geq 2$.} \\ [1ex]

   \end{array}
\right.
$$

\end{theorem}
\begin{proof} We consider the two cases separately.
\begin{enumerate}
\item $q=k(a-1)$ with $k\geq 1$. Then $ a \le q+1$. From Lemma~\ref{lem:broadcast-at-most-qa}, $\overrightarrow{C}(qa;1,a))$ admits an $(a-1)$-bounded independent broadcast. Moreover, since $n=ka(a-1)$ and from Theorem~\ref{th:ex 1act 1} we get
 $$\beta_b(\overrightarrow{C}(qa;1,a)) = \beta_b(\overrightarrow{C}(qa;1,-(q-1)a))  = ka(a-2)=(q + k)(a-2).$$
\item If $a=kq+1$ with $k\geq 2$, then $q+1 < a$. From Lemma~\ref{lem:broadcast-at-most-qa}, $\overrightarrow{C}(n;1,a))$ admits a $q$-bounded independent broadcast. Now, from Proposition~\ref{prop:bound-Inf-qa}, we have
$$\beta_b(\overrightarrow{C}(q(kq+1);1,a)) \geq a(q-1).$$

Moreover, from Proposition~\ref{prop:bound2-q} we get that

$$\sigma(f) = f(V) \leq \left\lfloor  \left( \frac{q-1}{q} \right) \left(n - \left|V_f^{1}\right| \right) \right\rfloor - (q-2)\left|V_f^{1}\right|,$$
for every $q$-bounded independent broadcast $f$ on $\overrightarrow{C}(n;1,a)$,
which implies $$\beta_b(\overrightarrow{C}(n;1,a)) \le a(q-1).$$

Finally, by Observation~\ref{isomorphisme}, we then get
$$\beta_b(\overrightarrow{C}(qa;1,a)) = \beta_b(\overrightarrow{C}(qa;1,-(q-1)a)) = a(q-1).$$
\end{enumerate}
This completes the proof.
\end{proof}

\begin{theorem}
If $n$, $a$, $q$ and $k$ are four integers such that $n=qa + a-1 $,  $3 \le q$ and $4 \le a$, then

$$ \beta_b(\overrightarrow{C}(n;1,a))  = \beta_b(\overrightarrow{C}(n;1,-(qa-1))  =
\left\{
   \begin{array}{ll}
     (q + k + 1)(a-2) & \text{if $q=k(a-1)$, $ k\geq 1$, } \\ [1ex]
     (a + k)(q-1) & \text{if $a=kq+1$, $k\geq 2$.} \\ [1ex]
   \end{array}
\right.
$$

\end{theorem}
\begin{proof}
Let $n=qa+a-1$. We consider the two cases separately.
\begin{enumerate}
\item If $q=k(a-1)$ with $k\geq 1$, then $ a \le q+1$. From Lemma~\ref{lem:broadcast-at-most-q(a-1)-1}, $\overrightarrow{C}(n;1,a))$ admits an $(a-1)$-bounded independent broadcast. Moreover, since $n=(ka+1)(a-1)$ and from Theorem~\ref{th:ex 1act 1} we get
 $$\beta_b(\overrightarrow{C}(n;1,a))  = (ka+1)(a-2)=(q + k + 1)(a-2).$$

\item If $a=kq+1$ with $k\geq 2$, then $q+1 < a$. From Lemma~\ref{lem:broadcast-at-most-q(a-1)-1}, $\overrightarrow{C}(q(n;1,a))$ admits a $q$-bounded independent broadcast. Now, from Proposition~\ref{prop:bound-Inf-qa-1}, we have
$$\beta_b(\overrightarrow{C}(n;1,a)) \geq (a+k)(q-1)$$

Moreover, from Proposition~\ref{prop:bound2-q}, we get that

$$\sigma(f) = f(V) \leq \left\lfloor  \left( \frac{q-1}{q} \right) \left(n - \left|V_f^{1}\right| \right) \right\rfloor - (q-2)\left|V_f^{1}\right|,$$
for every $q$-bounded independent broadcast $f$ on $\overrightarrow{C}(q(n;1,a)$,
which implies $$\beta_b(\overrightarrow{C}(n;1,a)) \le (a+k)(q-1).$$

Finally, by Observation~\ref{isomorphisme}, we then get
$$\beta_b(\overrightarrow{C}(q(n;1,a)) = \beta_b(\overrightarrow{C}(q(n;1,a-n)) =  (a+k)(q-1).$$
\end{enumerate}
\end{proof}

\end{document}